\documentclass[smallextended]{svjour3} 

\RequirePackage{fix-cm}

\usepackage{epsfig}
\usepackage{epstopdf}
\usepackage{amsmath}
\usepackage{amsfonts}
\usepackage{bm}
\usepackage{todonotes}
\newcommand{\be}{\begin{equation}}
\newcommand{\ee}{\end{equation}}

\begin{document}

\title{Blow up in a periodic semilinear heat equation}

\author{Marco Fasondini \and John King \and J.A.C.~Weideman}


\institute{Marco Fasondini \at School of Computing and Mathematical Sciences\\ University of Leicester \\  LE1 7RH Leicester \\ United Kingdom \\
\email{m.fasondini@leicester.ac.uk} \and  J.R. King \at School of Mathematical Sciences and Synthetic Biology Research Centre\\ University of Nottingham\\  NG7 2RD Nottingham\\ United Kingdom \\ \email{John.King@nottingham.ac.uk} \and
	J.A.C.~Weideman \at
	Department of Mathematical Sciences \\
	Stellenbosch University\\
	Stellenbosch 7600 \\
	South Africa\\
	\email{weideman@sun.ac.za}  
}


%


\maketitle

\abstract{Blow up in a one-dimensional semilinear heat equation is 
studied using a 
combination of numerical and analytical tools. The focus is on problems periodic
in the space variable and starting out from a nearly flat, positive initial condition.
Novel results include various asymptotic approximations that are, in combination,
valid over the entire
space and time interval right up to and including the blow-up time. 
Preliminary results on continuing a numerical solution beyond the singularity are also presented.}
 



\section{Introduction}\label{sec:intro}
The blow-up phenomenon in nonlinear evolution equations has been studied extensively
in the literature.   Some studies
focus on physical applications, such as singularity formation in 
fluids~\cite{braun2,braun1,hocking,Lushnikov21}, runaway in thermal processes~\cite{Dold91,Herrero,Lacey83},
and biological applications~\cite{jabbari2013discrete}.  Others 
deal with numerical aspects such as the computation
of blow-up profiles, estimation of blow-up times, and singularity 
tracking~\cite{Berger88,Budd96,Keller1993,Tourigny94,W03}.
For a general review, see~\cite{Galaktionov2002}.

The present paper is a continuation of the papers~\cite{Keller1993} and \cite{W03}.
The equation considered in these papers is the  
nonlinear heat equation
\begin{equation}
u_t = u_{xx} + u^2
\label{eq:u}
\end{equation}
(although more general nonlinearities and more than one space dimension were also considered in~\cite{Keller1993}).   In this paper we consider the initial condition
\begin{equation}
u(x,0)  = \frac{1}{\alpha-\epsilon \cos x}, \quad 
0 < \epsilon \ll \alpha,  
\label{eq:ic}
\end{equation}
i.e.~periodic, positive, and nearly flat.  (A broader class of nearly flat
initial data is considered in Appendix~\ref{sec:matchedasympt}.)
Figure~\ref{fig:schematic} shows a typical blow-up scenario
for this case.

\begin{figure}[htb]  
	\begin{center}
		\includegraphics[width=0.6\textwidth]{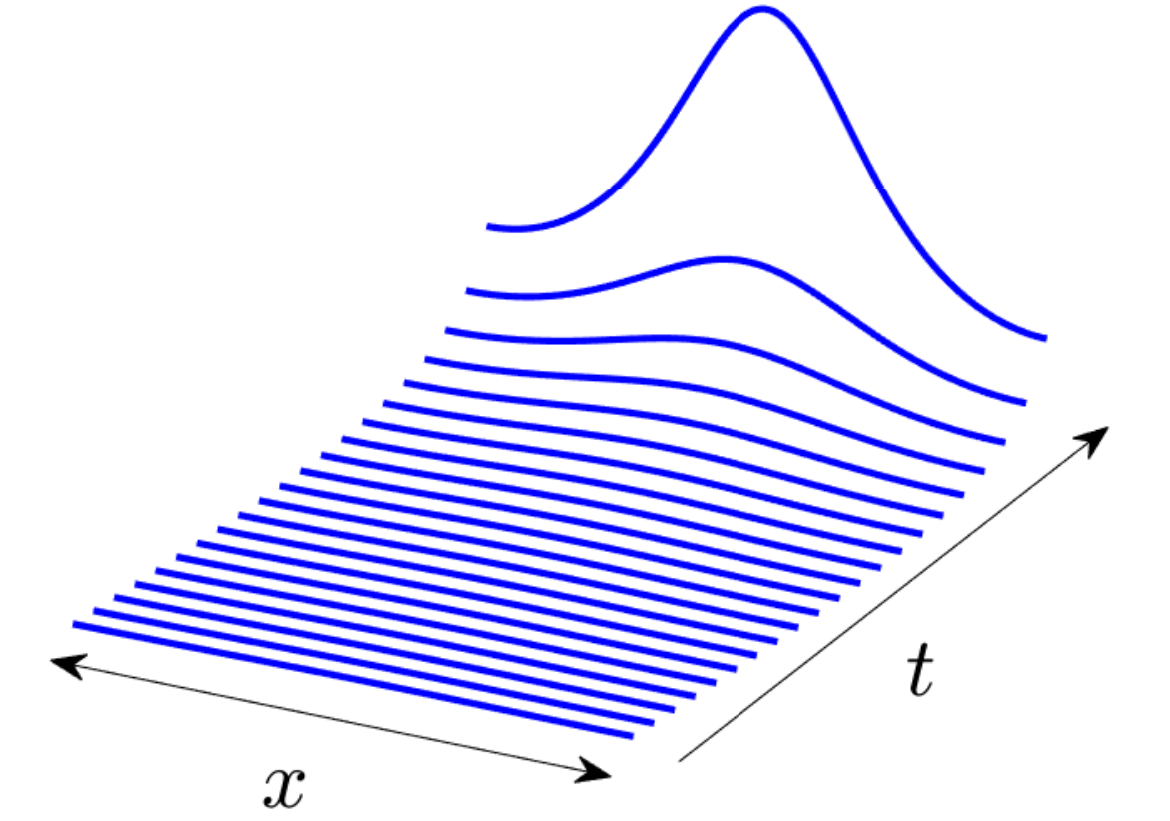}
	\end{center}
	\caption{Typical solution profiles of 
		equation~(\ref{eq:u}), in the case of a nearly flat initial condition such as~(\ref{eq:ic}). 
		The solution stays almost flat
		for a long while until a point blow-up occurs over a relatively short period of time.  }  
	\label{fig:schematic}
\end{figure}

In~\cite{W03}, it was shown by numerical computation that the
approach to blow up is not necessarily uniform.   That is, there may be times when the diffusive term dominates, leading to a
flattening of the solution profile.  At other times, particularly
near blow up, 
the nonlinearity dominates, leading to  a
steepening of the profile.   By numerically continuing the solution into  the complex plane,  it was shown that this 
behaviour can be associated with the dynamics of the
complex singularities of the solution.  When diffusion
dominates, they typically move farther from the real axis while
the opposite is true when nonlinearity dominates.   A 
point blow up occurs when the singularities reach the
real axis.  For the nearly flat initial conditions considered in this paper, however we show 
 that there is no  simple correspondence between the flatness of the solution profile on the real interval $[-\pi, \pi]$ and the proximity of the nearest singularity, except in the blow-up limit, where the steepening solution profile is associated with the impingement of singularities onto the real axis at $x=0$.

In~\cite{Keller1993}, the focus was not on complex singularities but on 
asymptotic estimation of the blow-up time and the solution profile
near blow up.  The starting point in that paper was
the substitution $u = 1/v$, which transforms~(\ref{eq:u}) into
\begin{equation}
v_t = v_{xx}-1-2(v_x)^2/v.
\label{eq:v}
\end{equation}
This transformation has advantages for both analysis and computation.   
Advantages for analysis are spelled out
in~\cite{Dold91}.   For numerical computation, it is easier
to deal with solutions tending to zero than to infinity.  
This avoids the need for   specialized rescaling
algorithms or moving mesh methods~\cite{Berger88,Budd96}.
The downside is that the simple polynomial
nonlinearity in~(\ref{eq:u}) has been replaced by the 
more complicated nonlinearity in~(\ref{eq:v}). While this can be ameliorated by multiplication by $v$, it
raises the red flag of division by zero as $v$ becomes small.
However, it was shown in~\cite{Keller1993} and reaffirmed in section~\ref{sec:blowuplim} of the
present paper that the 
nonlinear term  remains bounded as $v$ approaches
zero.  This further suggests the possibility of integration through
the zero of $v$, i.e.~through the singularity of $u$, but~\cite{Keller1993} reported a failed effort.  We continued those investigations and announce  preliminary
findings here.

In~\cite{W03} periodic boundary conditions were considered
while~\cite{Keller1993} looked at the pure initial value
problem on the entire real line. In this paper we continue with the
periodic case, as this gives us access to highly accurate
Fourier spectral methods (which can also be applied to 
problems on the infinite line, but with a substantial penalty
in accuracy).  One contribution here is a conversion
of the analysis of~\cite{Keller1993} to the periodic situation, 
which is not just a triviality but contributes significant
new results as we shall discuss.

For the numerical computations of this paper,   a full spectral method
based on a Fourier series
\begin{equation}
v(x,t) = \sum_{k=-\infty}^{\infty} c_k(t) e^{ikx}, \qquad -\pi \leq x < \pi,
\label{eq:fourier1}
\end{equation}
is used.   The derivatives on the right-hand side of~(\ref{eq:v}) 
are computed by analytical 
differentiation of this series, and the nonlinear terms by convolution
and de-convolution.  The result is an infinite system of ODEs for the evolution
of the Fourier coefficients $c_k$, which we truncate  
at $|k| = 128$.   To integrate this system, we use {\tt ode45},
\uppercase{MATLAB}'s standard ODE solver.  It has adaptive time stepping
to maintain accuracy, and the 
tolerance parameters for this are set to a stringent $10^{-12}$. To compute
the time $t = t_c$ at which the numerical solution blows up 
(i.e.~$v=0$), we use the fact that the initial conditions considered here lead to blow up at $x=0$ and hence we check when
the sum of the $c_k$ equals zero.  This can be done by the `event' option in  {\tt ode45}.  

We use this numerical solution as reference solution for the purposes of 
checking the various  asymptotic estimates.  
Because of the
relatively smooth nature of $v$, even at the critical time (as will be discussed), we believe it is sufficiently close to the true solution for all
verification purposes. 
 
The three main sections of the paper can be summarized in a nutshell as: 
before blow up, at blow up, and after blow up.   More specifically,
in section~\ref{sec:twomode} we
present a perturbation analysis that approximates
the solution to  (\ref{eq:v}) accurately to $O(\epsilon^2)$ on the entire
periodic space domain and the whole time interval until a time 
that is
$O(\epsilon)$  close to 
blow up.  Beyond that time up to the blow-up 
point it has to be modified
and this is done by 
matched asymptotic expansions, the details of which are contained in Appendix~\ref{sec:matchedasympt}. In this appendix we also analyse the dynamics of the singularities of the solution and in Appendix~\ref{sec:flatness} the relation between the proximity of the singularities to the real axis and the steepness of the solution profile on $[-\pi, \pi]$ is clarified. In sections~\ref{sec:twomode} and~\ref{sec:blowuplim}, the analyses in the appendices are confirmed by numerical experiments.
Section~\ref{sec:postblowup} is a preliminary report on 
our efforts in
integrating through the singularity at the blow-up time and the subsequent evolution.  

\section{Two-mode perturbation analysis}\label{sec:twomode}

The analysis of~\cite{Keller1993} was based on the truncated Taylor expansion 
\begin{equation}
v  \approx a(t)+b(t)x^2. \label{eq:Taylans}
\end{equation}
By substituting into~(\ref{eq:v})  
and dropping powers of $x^4$ the problem was reduced to a dynamical system in
the variables $a$ and $b$; see~(\ref{eq:twomodeKL}) below.   Here we follow an analogous procedure, but consider instead a truncated Fourier expansion 
\begin{equation}
v  \approx a(t)-b(t)\cos x.   
\label{eq:ansatz}
\end{equation}
Under the assumption of strictly 
positive solutions, i.e., $0 < b(t) < a(t)$,
both of these 
approximations blow up (in the variable $u = 1/v$) in finite time at $x=0$.

\begin{figure}[htb]  
	\begin{center}
		\includegraphics[width=0.95\textwidth]{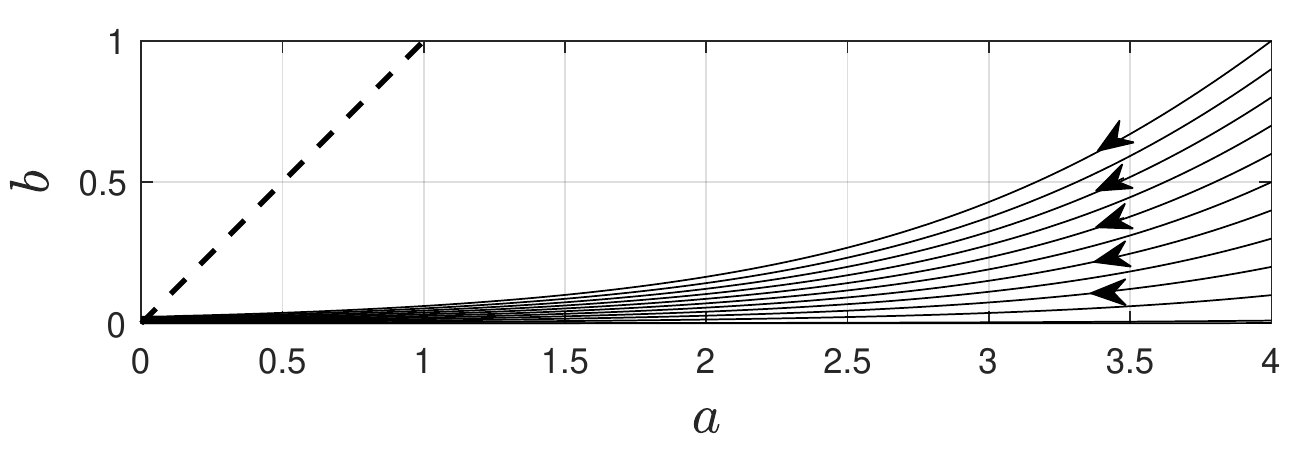}	
	\end{center}
	\caption{Phase plane of the system~(\ref{eq:twomode1}) in the domain
		$0 < b < a$.  The analysis of this section
	approximates solution curves close to $b=0$.  Blow up at $x=0$ in the 
	$u$-equation~(\ref{eq:u}) corresponds to solution trajectories
	intersecting the dashed line $b=a$. 
}  
	\label{fig:phase}
\end{figure}

Substitution of~(\ref{eq:ansatz}) into~(\ref{eq:v}) and neglecting
$\cos 2x$ contributions gives the system
\begin{equation}
a \frac{da}{dt} + \frac12 b \frac{db}{dt} = -a -\frac32 b^2, 
\quad b \frac{da}{dt} + a \frac{db}{dt} = -ab -b,
\label{eq:twomode}
\end{equation}
or, in explicit form,
\begin{equation}
\frac{da}{dt} = \frac{2ab^2+2a^2-b^2}{b^2-2a^2}, \quad 
\frac{db}{dt}  = \frac{b(2a^2-3b^2)}{b^2-2a^2}.
\label{eq:twomode1} 
\end{equation}
(Note that the assumption $0 < b < a$ precludes the vanishing of the denominators.)
The phase plane of this system is shown in~Figure~\ref{fig:phase}.  

\begin{figure}[htb]  
\mbox{		\includegraphics[width=0.495\textwidth]{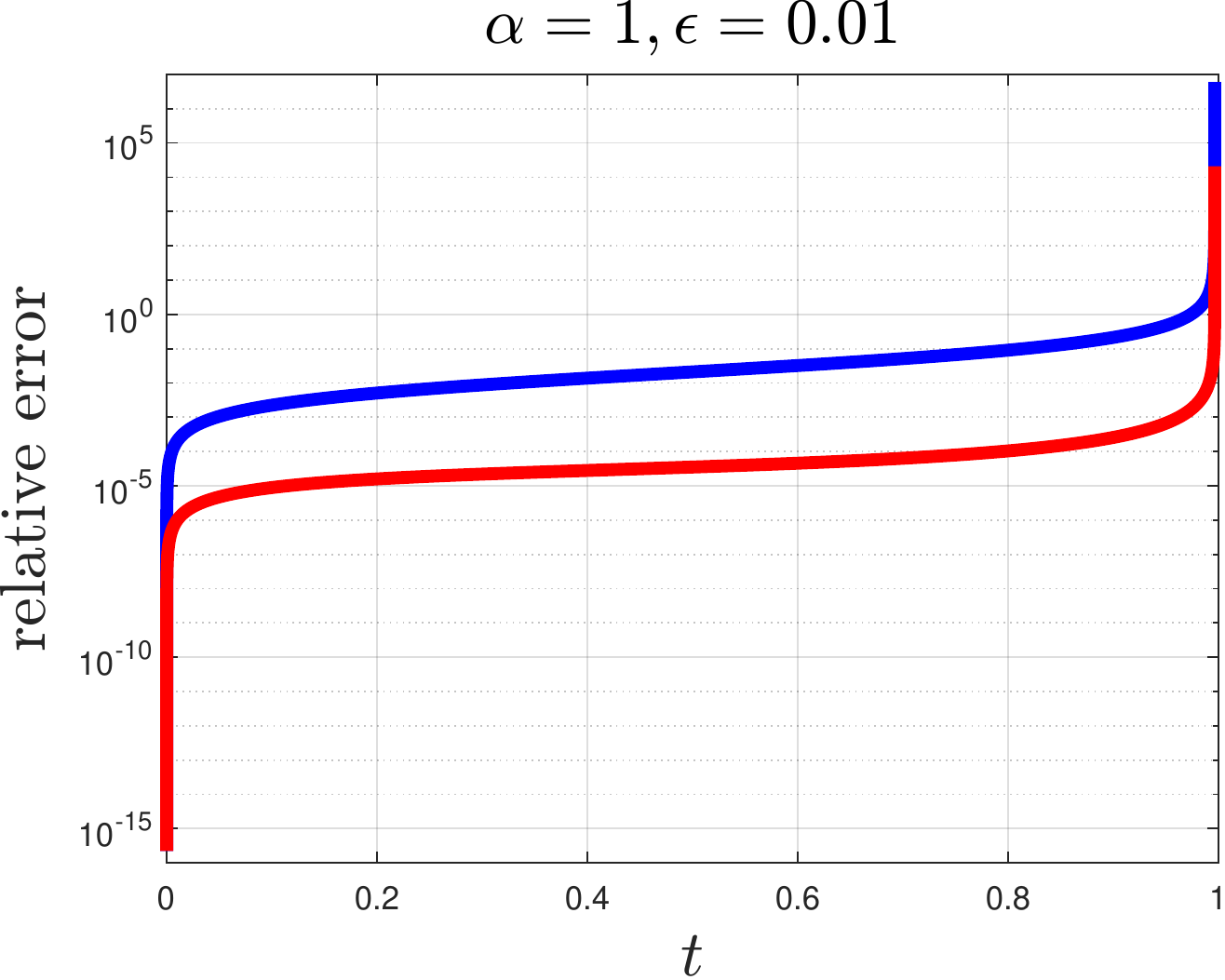}
		\includegraphics[width=0.495\textwidth]{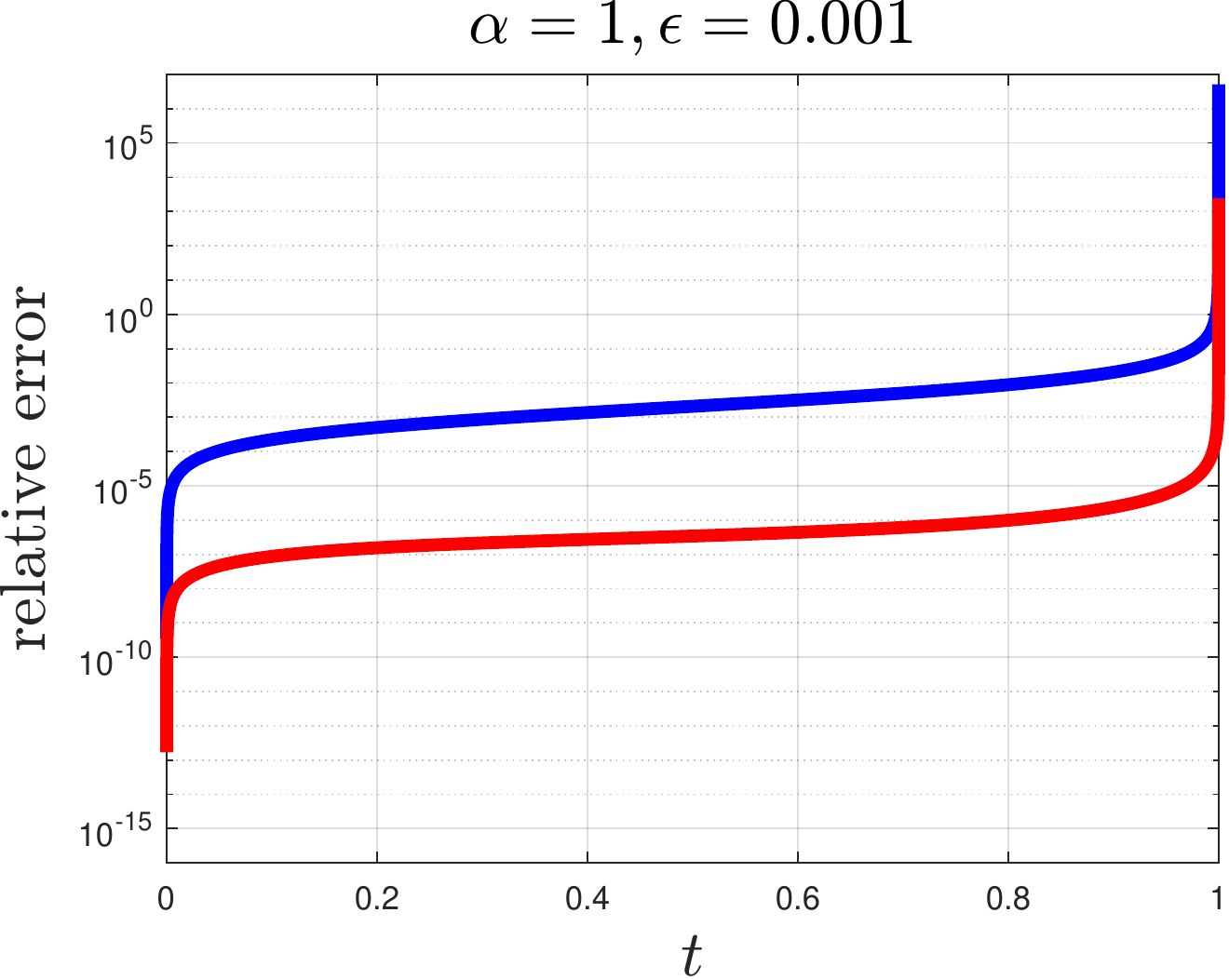}}
	\caption{The maximum relative error on $[-\pi, \pi]$ of the perturbation approximation (\ref{eq:pert}) (blue) and the numerical solution to the two-mode system (\ref{eq:twomode1}) (red)  as a function of $t$.  The errors are calculated with reference to the numerical solution of (\ref{eq:v})
	mentioned in the final paragraph of section~\ref{sec:intro}.
	} \label{fig:pertacc} 
\end{figure}

Consider solution curves in~Figure~\ref{fig:phase} that originate near $b = 0$, say
\begin{equation}
a(0) = \alpha, \quad b(0) = \epsilon, \quad 0 < \epsilon \ll \alpha.  \label{eq:FICs}
\end{equation}
As an explicit solution of the system~(\ref{eq:twomode}) appears not to exist, we
settle for a perturbation analysis, by expanding
\begin{equation}
a = a_0(t) + \epsilon a_1 (t) + 
O(\epsilon^2),   \qquad 
b = \epsilon b_1(t) + 
O(\epsilon^2) \label{eq:F1sttimescale}
\end{equation}
where $a_0(0) = \alpha$, $b_1(0) = 1$, $a_1(0)= 0$. 
Substitution into~(\ref{eq:twomode}) gives, to zeroth order
\begin{equation}
a_0 \frac{da_0}{dt} = -a_0 \quad \Rightarrow \quad a_0 = \alpha-t.  \label{eq:a0s}
\end{equation}
At $O(\epsilon)$,
\begin{equation}
\frac{da_1}{dt} = 0, \ \frac{db_1}{dt} = -b_1 \quad \Rightarrow \quad 
a_1 = 0, \ b_1  = e^{-t}.  \label{eq:b0s}
\end{equation}
This gives, to $O(\epsilon^2)$, the approximate solution
\begin{equation}
\widetilde{v} = \alpha -t - \epsilon \, e^{-t} \cos x.
\label{eq:pert}
\end{equation}
Plugging this expression into~(\ref{eq:v}) gives
\begin{equation}
\widetilde{v} \big( \widetilde{v}_t - \widetilde{v}_{xx}+ 1+ 2(\widetilde{v}_x)^2/\widetilde{v} \big)  = 2 \, \epsilon^2 e^{-2t} \sin^2 x,
\label{eq:remainder}
\end{equation}
which confirms that the $v$-equation is satisfied to $O(\epsilon^2)$ uniformly in $x$,
for all $\widetilde{v}$ bounded away from zero.   Figure~\ref{fig:pertacc}  illustrates that the perturbation approximation (\ref{eq:pert}) and the (numerically computed) solution to the two-mode system (\ref{eq:twomode1}) are good approximations to the numerical reference solution of (\ref{eq:v}) on most of the interval $[0, t_c)$, where $t_c$ is the critical time at which blow up occurs.
As $t \to t_c$, the assumption underlying the perturbation analysis and the two-mode approximation, $b(t) \ll a(t)$, is no longer valid, 
however, and the approximations lose accuracy. 


Various estimates can be obtained from the perturbation solution~(\ref{eq:pert}).
By setting $\widetilde{v}=0$ at $x=0$ and excluding $O(\epsilon^2)$ terms, for example, one obtains the following estimate for the blow-up 
time
\begin{equation}
t_c \approx \widehat{t}_c :=  \alpha -e^{-\alpha} \epsilon.
\label{eq:thatdef}
\end{equation}
In Appendix~\ref{sec:matchedasympt}, using the method of matched asymptotic expansions, a higher-order estimate of $t_c$ is derived (see (\ref{eq:beta12}) and (\ref{eq:Fourbetagamma})): 
\begin{equation}
    t_c \approx \widetilde{t}_c := 
    \alpha - e^{-\alpha}\epsilon -(2C_1+C_2+C_3)\epsilon^2, 
    \label{eq:ttildedef}
\end{equation}
the constants $C_i$ being defined in (\ref{eq:cdefs}).   The accuracy 
of these estimates is verified in Table~1. (The values of
$t_c$ listed in the table were computed by the method described in section~\ref{sec:intro}, and are believed to be correct to all digits shown.)


\begin{table}[htb]
\begin{center}
	\renewcommand{\arraystretch}{1.2}
 \begin{tabular}{c|cccc|}
 	& \multicolumn{4}{c|}{$\alpha = 0.25$}   \\
$\epsilon$ & 	$t_c$ & 	
$t_c'-t_c$  & $\widehat{t}_c-t_c$ & $\widetilde{t}_c - t_c$   \\  \hline 
0.1 & 0.161963 &  $-$3.6e-04 & 1.0e-02 & 2.6e-02  \\
0.01 & 0.242093  &  1.8e-03 &  1.2e-04  &  2.8e-04  \\
0.001 & 0.249220 & 2.2e-04  & 1.2e-06  & 2.8e-06  \\
\hline
	& \multicolumn{4}{c|}{$\alpha = 1$}   \\
	\hline
0.1 & 0.955542 & 2.1e-03 &  7.7e-03 & 4.5e-03  \\
0.01 & 0.996241 & 9.5e-04 & 8.1e-05 &  4.9e-05   \\
0.001 & 0.999631 & 1.1e-04 & 8.1e-07  & 4.9e-07 \\
\hline
 	& \multicolumn{4}{c|}{$\alpha = 4$}   \\
 	\hline
0.1 & 3.996685 &  5.0e-04 & 1.5e-03  &  1.4e-05   \\
0.01 & 3.999802  & 5.3e-05 & 1.5e-05  &  1.2e-07  \\
0.001 & 3.999982 & 5.4e-06 & 1.5e-07  &  1.2e-09  	
 \end{tabular}
\end{center}
\caption{Blow-up times for various parameter choices in the initial
	condition~(\ref{eq:ic}):   $t_c$ is the blow-up time as computed
from the reference solution, $t_c'$ is the blow-up time
as estimated from a numerical solution of the two-mode system~(\ref{eq:twomode}) and $\widehat{t}_c$ and $\widetilde{t}_c$ denote the 
estimates (\ref{eq:thatdef}) and (\ref{eq:ttildedef}), respectively.
}
\label{tab:tc}
\end{table}

\begin{figure}[htb]  
	\mbox{	
	\includegraphics[width=0.495\textwidth]{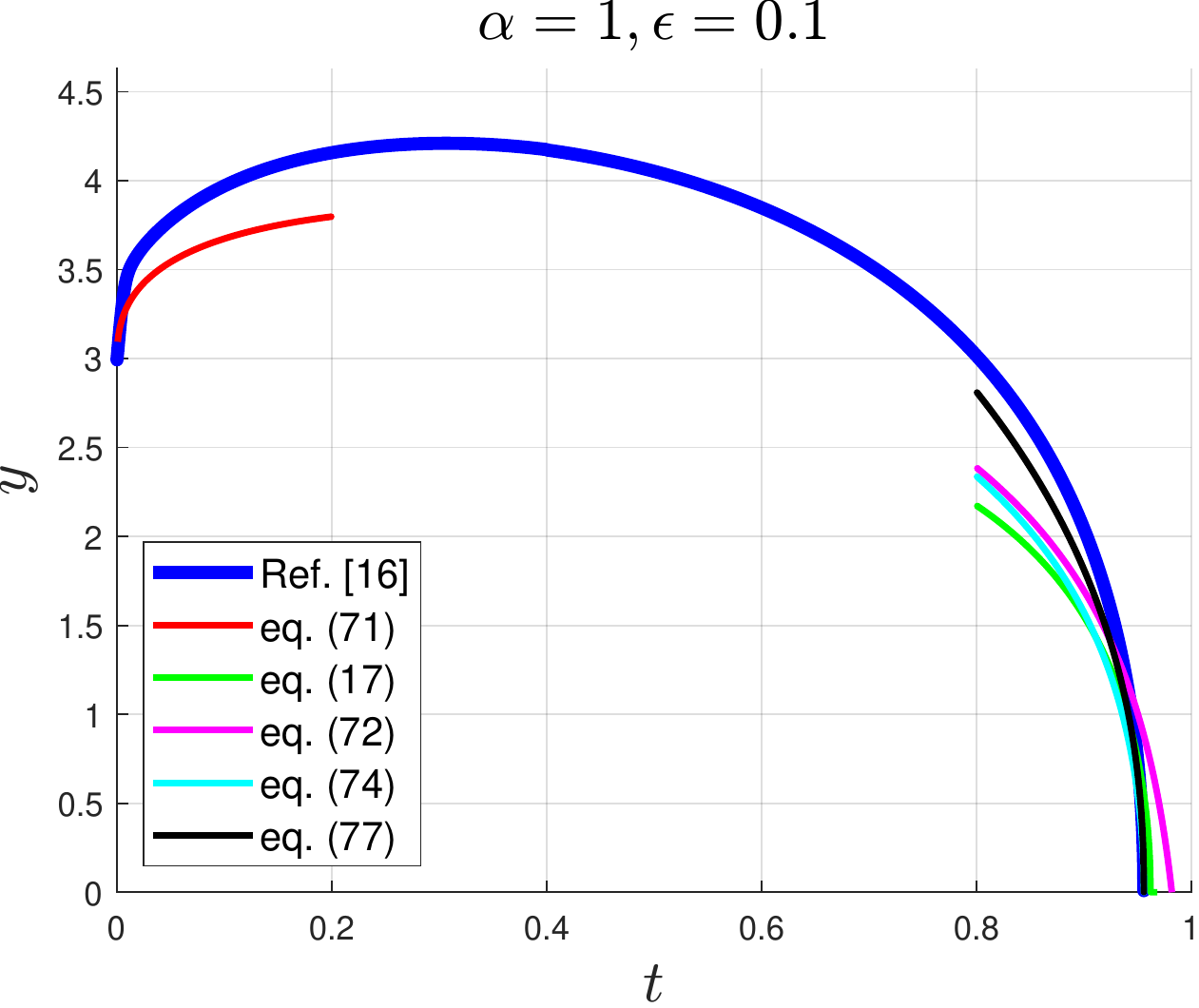}
	\includegraphics[width=0.495\textwidth]{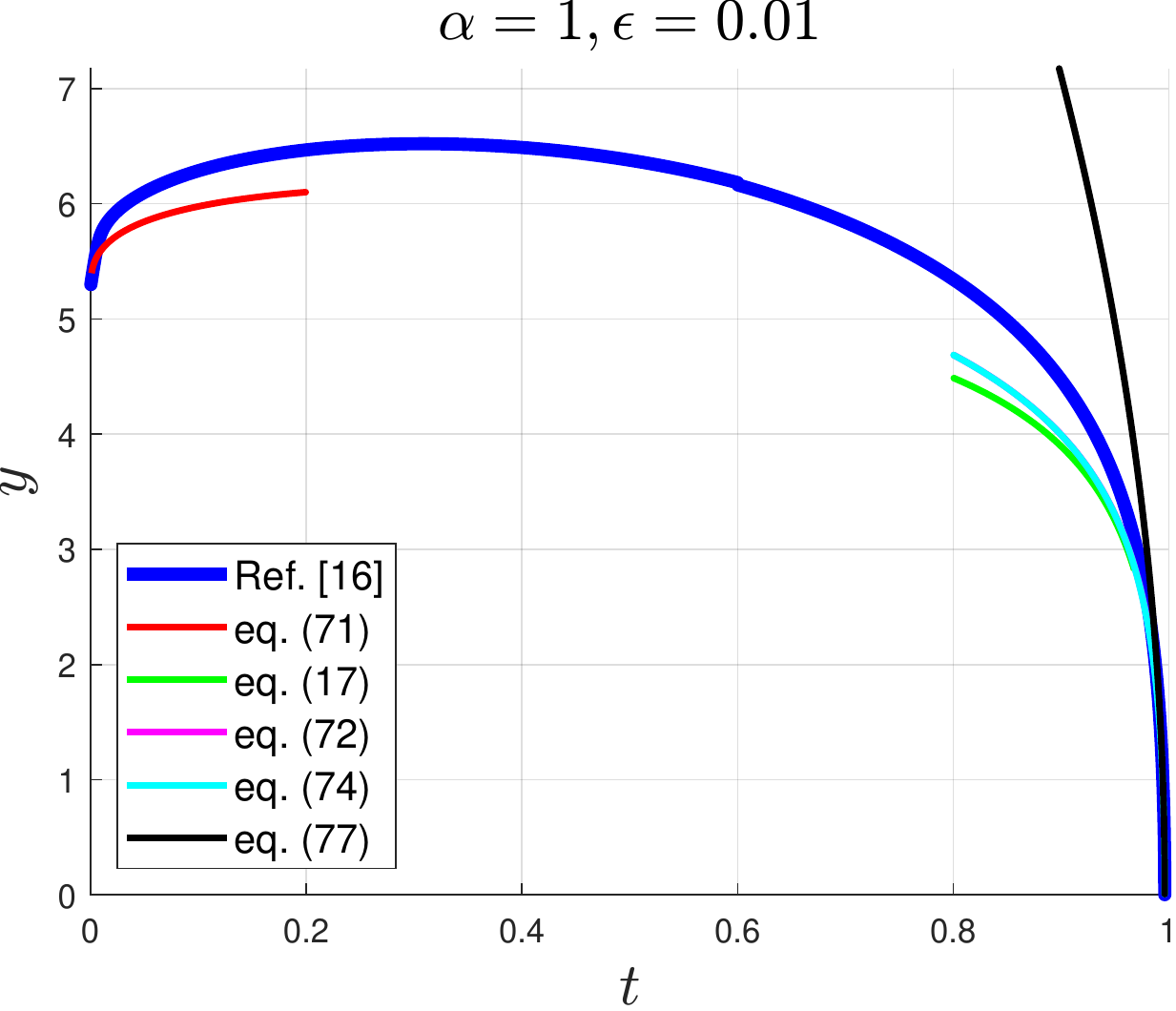}}
	\caption{The position on the positive imaginary axis of the singularity of the solution to (\ref{eq:v}) with initial condition $v(x,0) = \alpha - \epsilon\cos x$. (In the right frame, the approximation (\ref{eq:y1stscalel2}) is not visible because it is indistinguishable from (\ref{eq:y2ndscale}).)
}  
	\label{fig:YvsT}
\end{figure}

The singularity dynamics mentioned in the introduction can be estimated
as follows. The complex singularity in the $u$-equation 
corresponds to a complex zero of the $v$-equation.  Assuming it is located at
$x = iy$ with $y$ real, setting $\widetilde{v}=0$  gives
\begin{equation}
y = \cosh^{-1} \big((\alpha-t) e^{t}/\epsilon\big).
\label{eq:YvsT}
\end{equation}
This approximation is valid, however, only for small values of $y$, i.e., near blow up.  This follows from the fact that the remainder term on the right-hand side 
of~(\ref{eq:remainder}) 
grows exponentially with $y$. In Figure~\ref{fig:YvsT}, the asymptotic estimate (\ref{eq:YvsT}) is compared to a numerical estimate of the  singularity location obtained via the method  of~\cite{sulem1983}, which is essentially an estimation of the width of the strip of analyticity of the solution, by examining the rate of decay of its Fourier
coefficients. To apply the method of~\cite{sulem1983}, we use the fact that, to leading order and away from the blow-up time, singularities of solutions to the $u$-equation are second-order poles\footnote{In~\cite{FKW} we show that these singularities are in fact logarithmic branch points, however the branch point singularity appears in the fourth-order term in the local expansion about the singularity. Therefore, for the purpose of estimating the position of the singularity, we only use its leading-order, second-order pole behaviour. We shall find in section~\ref{sec:blowuplim} that in the limit $t \to t_c$, the leading order behaviour of the singularity at $x=0$, which results from the coalescence of two singularities, is of a more complicated form than that of a second-order pole.   }. 
Figure~\ref{fig:YvsT} also shows asymptotic estimates of the singularity location that are derived in Appendix~\ref{sec:comp} by the method of matched asymptotic expansions. The estimates are valid in the limits $t \to 0^{+}$ and $t \to t_c^{-}$ but in between, for $t = O(1)$, the singularity location of the original PDE (\ref{eq:v}) is described by the singularity location of a more difficult nonlinear initial value problem (\ref{eq:backdiff})--(\ref{eq:backdiffcond2}) which is not analytically solvable. The asymptotics nonetheless correctly indicate the initial movement of the singularity away from the real axis (at a speed that becomes infinite as $t\to 0^{+}$) and the final motion shortly before the singularity collides with the real axis at $t = t_c$.

In Appendix~\ref{sec:flatness} we investigate the relationship between the position of the complex singularity and the height of the peak of the solution profile in the $u$-variable (which is located at $x = 0$, see Figure~\ref{fig:schematic}) relative to the solution value at $x = \pm \pi$.




Returning to the analysis of~\cite{Keller1993}, which was based on the truncated Taylor 
approximation 
(\ref{eq:Taylans}), we note that the focus in that paper was on the behaviour
at $t=t_c$, not the evolution on $[0,t_c)$ as is the focus
here.  It is therefore instructive to adapt that analysis here and 
compare results to~(\ref{eq:pert}) and (\ref{eq:remainder}). 

The dynamical system analogous to~(\ref{eq:twomode})--(\ref{eq:twomode1}) is now
\begin{equation}
\frac{da}{dt} = 2b-1, \quad a\frac{db}{dt}  = -8 b^2.
\label{eq:twomodeKL}
\end{equation}
This is a much simpler system, and in fact admits a first integral 
$2 \log b + 1/b + b \log a = \mbox{constant}$, although we shall not make
use of this result and neither was it used in~\cite{Keller1993}.  (No
such first integral could be found for~(\ref{eq:twomode})--(\ref{eq:twomode1}).)

Proceeding with a perturbation analysis based on~(\ref{eq:FICs})--(\ref{eq:F1sttimescale}) give
$a_0 = \alpha-t$, $a_1(t) = 2$, and $b_1(t) = 1$. Therefore, excluding $O(\epsilon^2)$ terms, $v$ is approximately
\begin{equation}
\widetilde{v} = \alpha-t + \epsilon (2t+x^2)   
\label{eq:pertKL}
\end{equation}
in analogy with~(\ref{eq:pert}).  Plugging into~(\ref{eq:v}) gives
\begin{equation}
\widetilde{v} \big( \widetilde{v}_t - \widetilde{v}_{xx}+ 1+ 2(\widetilde{v}_x)^2/\widetilde{v} \big)  =  8 \, \epsilon^2 x^2.
\label{eq:remainderKL}
\end{equation}
Comparing with~(\ref{eq:remainder}) shows an advantage of the periodic
analysis, namely that its right-hand side is $O(\epsilon^2)$ uniformly
in $x$, whereas the right-hand side of~(\ref{eq:remainderKL}) is
 $O(\epsilon^2)$ only if $x^2 = O(1)$.  That is, with $\widetilde{v}$ bounded away from zero
 (\ref{eq:pert}) provides a valid approximation over
 the entire (periodic) domain, while~(\ref{eq:pertKL}) is valid
 only near $x=0$.


 On the other hand, the advantage of the analysis of~\cite{Keller1993}
 is that it gives a valid description of the structure
 of the solution close to blow up, as discussed in  
  Appendix~\ref{sec:matchedasympt}.

\section{Solution in the blow-up limit }\label{sec:blowuplim}

Appendix~\ref{sec:matchedasympt} is devoted to an asymptotic analysis of the solution on three time scales that are progressively closer to the blow-up time, namely $t = O(1)$, $T = (t-t_c)/\epsilon = O(1)$ and  $\tau = -\epsilon \log(-T) = -\epsilon\log((t_c-t)/\epsilon) = O(1)$. The analysis is performed using the method of matched asymptotics and the results on the first time scale are the same as the those obtained in section~\ref{sec:twomode} via a regular perturbation analysis (in particular, recall the approximation (\ref{eq:pert}) and its failure close to the blow up time, as seen in Figure~\ref{fig:pertacc}).

\begin{figure}[h!]  
\centering
\includegraphics[width=0.6\textwidth]{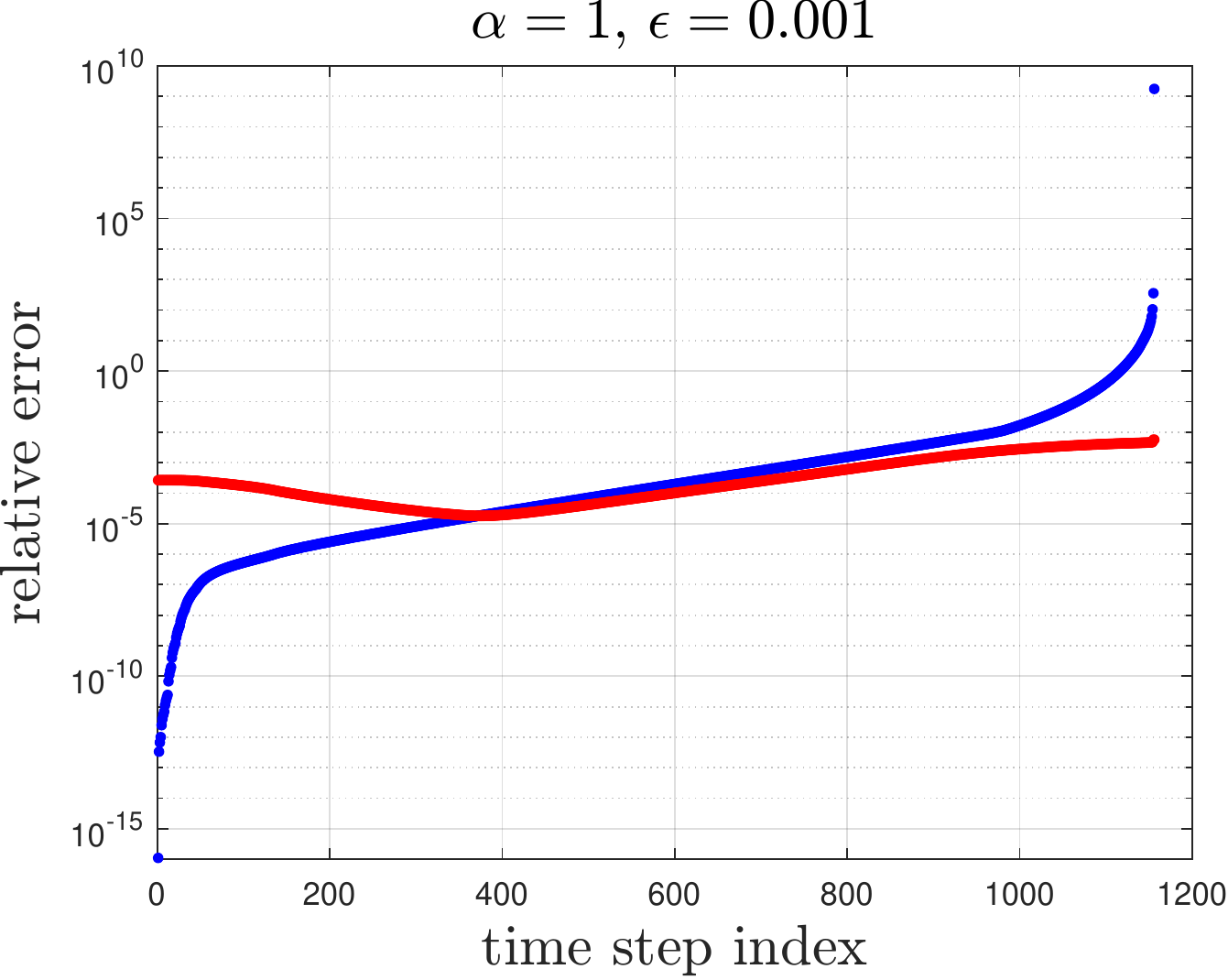}
	\caption{Errors in  the asymptotic approximations (\ref{eq:pert}) (blue) and (\ref{eq:timescale2}) (red) when compared to the reference solution. The numerical integration method returns approximations 
	at times $t = t_k$, $k = 0, 1, \ldots, M$, where $t_0 = 0$ and $t_M =t_c$.  Here $M = 1155$ and $t_c = 0.9996$ (more digits are listed
	in Table~\ref{tab:tc}).
	At  each $t = t_k$, the relative errors of (\ref{eq:pert}) and (\ref{eq:timescale2}) are calculated on the   interval $x \in [-\pi, \pi]$ and plotted against the time step index $k$. Because of the adaptive time-stepping, the
	$t_k$ values are not equidistant, but the spacing is much denser near 
	$t=t_c$. As an indication of this, note that the two curves intersect at $k = 372$ and $t_{372} \approx 0.973$,  which is already quite close to the critical time $t_c$ even though roughly 800 more time steps are to be taken. From the asymptotic analysis, we expect the approximation on the first time scale (\ref{eq:pert}) (blue) to start to break down and (\ref{eq:timescale2}) (red) to be valid on the second time scale when $(t_c - t)/\epsilon = O(1)$. Here we have that at around $k = 800$, $t_c - t_k \approx \epsilon = 0.001$, which is indeed when the error curves of the two approximations start to diverge noticeably.
	} \label{fig:allasympt}
\end{figure}

\begin{figure}[h!]  
	\mbox{
	\includegraphics[width=0.485\textwidth]{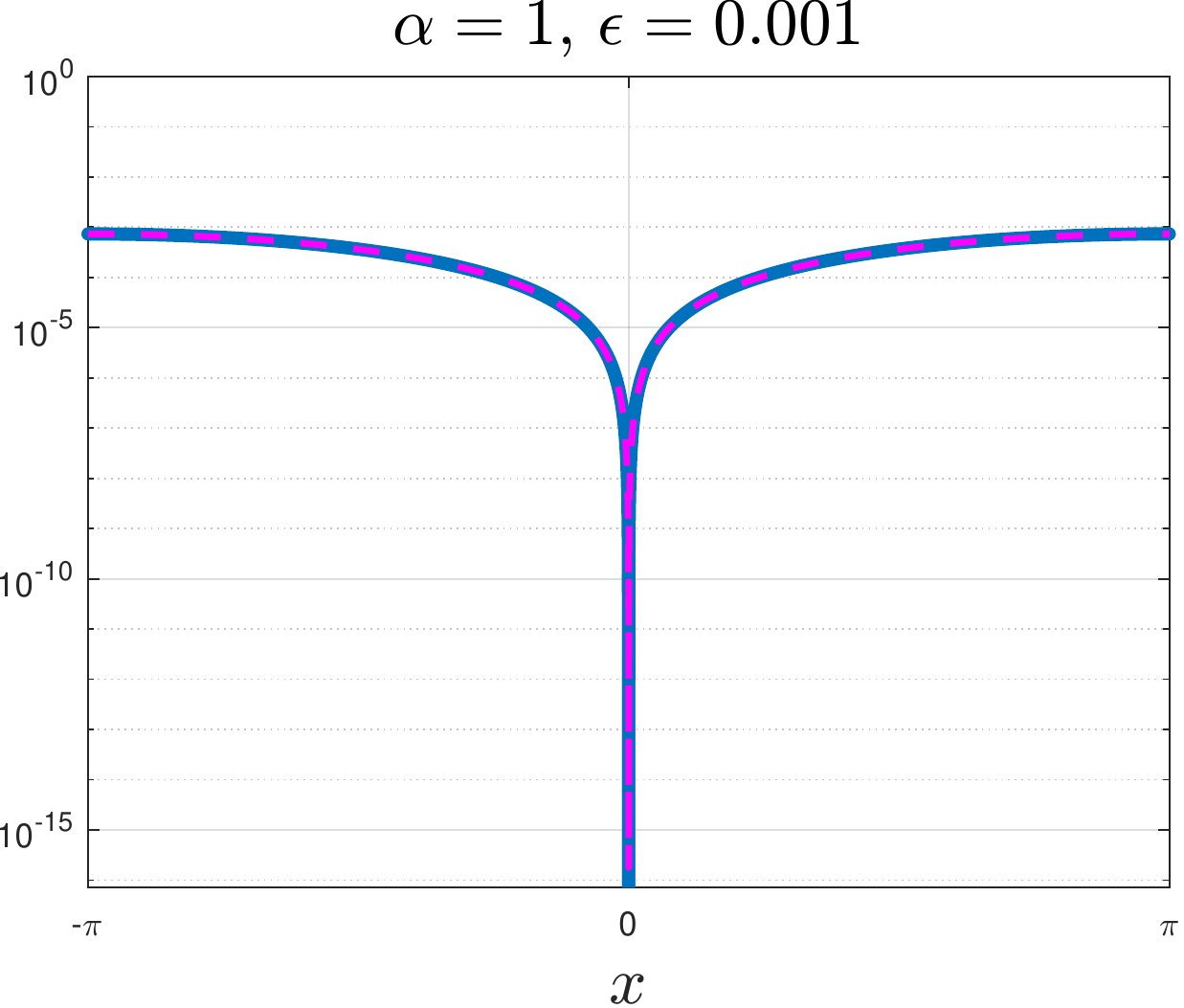}
		\includegraphics[width=0.51\textwidth]{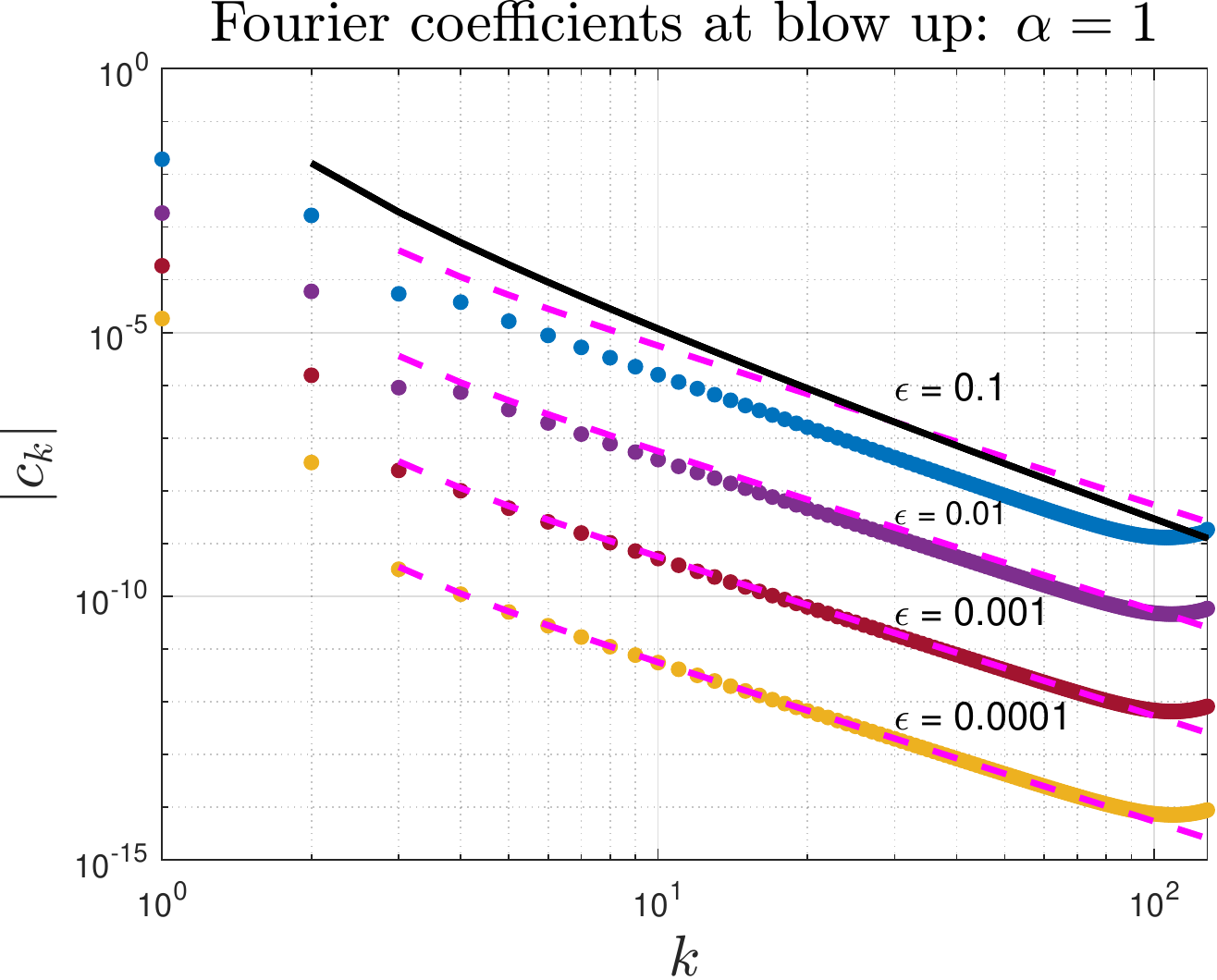}
	}
	\caption{
	Left: The solution to the $v$-equation (\ref{eq:v}), displayed on
	a log-scale, 
	at the blow-up time $t_c$ (solid curve) compared to the asymptotic estimate (\ref{eq:blowupprofsimp}) 
	(dashed curve on top of the solid curve). 
	Here we would like to draw attention to the fact that, at the origin
	and at the critical time, the values of $v$ are computed to 
	levels close to roundoff error ($\sim 10^{-16}$), 
	which translates into huge values in the $u$ solution near blow up.  
	(See for example the third frame in Figure~\ref{fig:PostB}.)
	This is a consequence of solving equation~(\ref{eq:v}) instead 
	of~(\ref{eq:u}).  
	Right: The Fourier coefficients of the
	reference numerical solution at the blow-up time (dots) compared to the coefficients of the global asymptotic blow-up profile (dashed lines),   given in (\ref{eq:BUFourcoeffs}). The single solid line shows the asymptotic estimate of the decay of the Fourier coefficients of the local blow-up profile, which is given in (\ref{eq:FourBUc}).
	} \label{fig:blowupprofile}
\end{figure}

On the second time scale, $T = O(1)$, the following asymptotic approximation is derived 
for $x = O(1)$: 
\begin{equation}
\begin{split}
& v \sim t_c-t + 2\epsilon\, e^{-\alpha}\sin^2(x/2)  +2\epsilon^2 \log \epsilon\, e^{-2\alpha}\sin^2 x +  \epsilon (t - t_c)e^{-\alpha}\cos x  \\
& + 2\epsilon^2\sin^2 x\bigg(e^{-2\alpha} \log\left( \frac{t_c-t}{\epsilon} + 2e^{-\alpha}\sin^2(x/2)   \right) + C_1 + C_3\bigg),
\end{split}\label{eq:timescale2}
\end{equation}
which is obtained by combining (\ref{eq:w0s})--(\ref{eq:w2s}) and (\ref{eq:phi12Four})--(\ref{eq:Fourbetagamma}). It also follows that, by setting $t = t_c$ in (\ref{eq:timescale2}), we obtain a representation of the solution profile at the blow-up time that is valid as $t \to t_c^{-}$, $\epsilon \to 0$ with $x = O(1)$:
\begin{equation}
v \sim  2 \epsilon\, e^{-\alpha} \sin^2 (x/2) +   
    2 \epsilon^2 \sin^2 x  \Big(e^{-2\alpha}   
   \log\left(2 \epsilon\, e^{-\alpha}\sin^2(x/2)\right)    + C_1 + C_3 \Big).
\label{eq:blowupprofsimp}
\end{equation}
Figure~\ref{fig:allasympt} verifies the accuracy of the asymptotic approximations (\ref{eq:pert}) (away from the blow-up time) and (\ref{eq:timescale2})--(\ref{eq:blowupprofsimp}) (close to and at the blow-up time).  The left frame of Figure~\ref{fig:blowupprofile} shows the numerical solution at the blow-up time with the blow-up profile (\ref{eq:blowupprofsimp}) superimposed on it. Hence, we  have accurate asymptotic expressions for the solution on the entire spatial interval $x \in [-\pi, \pi]$ and from $t = 0$ all the way up to and including the blow-up time. For $x \to 0$ at $t = t_c$, however, we shall need another asymptotic approximation, namely (\ref{eq:BUexpsmallx}), to be discussed below.

The rate of decay of the   Fourier coefficients  of  $v$    will be of relevance in the next section when we consider the possibility of continuing the solution beyond blow-up. From (\ref{eq:blowupprofsimp}) we deduce that the  $c_k$  decay as
$O(k^{-3})$ close to the blow-up time because
\[
\frac{1}{2\pi}\int_{-\pi}^{\pi} \sin^2\!x\,\log\left(\sin^2(x/2)\right) e^{-ikx} d x = \frac{2}{k(k^2 - 4)}
\]
which implies that
\begin{equation}
 c_k(t_c) \sim  
 \frac{4\epsilon^2 e^{-2\alpha}}{k^3}, \qquad k \to \infty. \label{eq:BUFourcoeffs}
\end{equation}
The right frame of Figure~\ref{fig:blowupprofile} confirms the accuracy of this estimate as $\epsilon \to 0$.



\begin{figure}[h!]  
	\mbox{
	\includegraphics[width=0.495\textwidth]{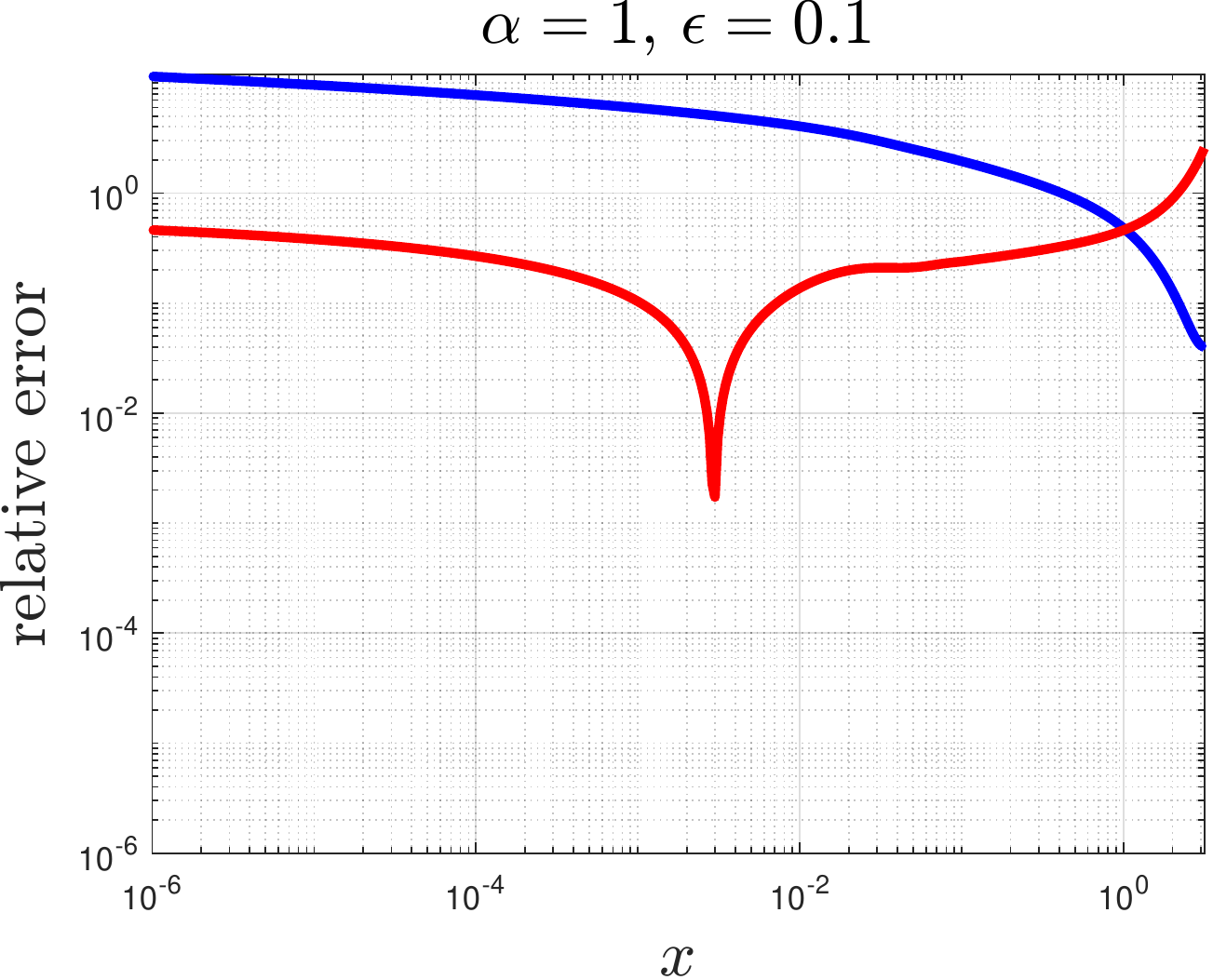}
		\includegraphics[width=0.495\textwidth]{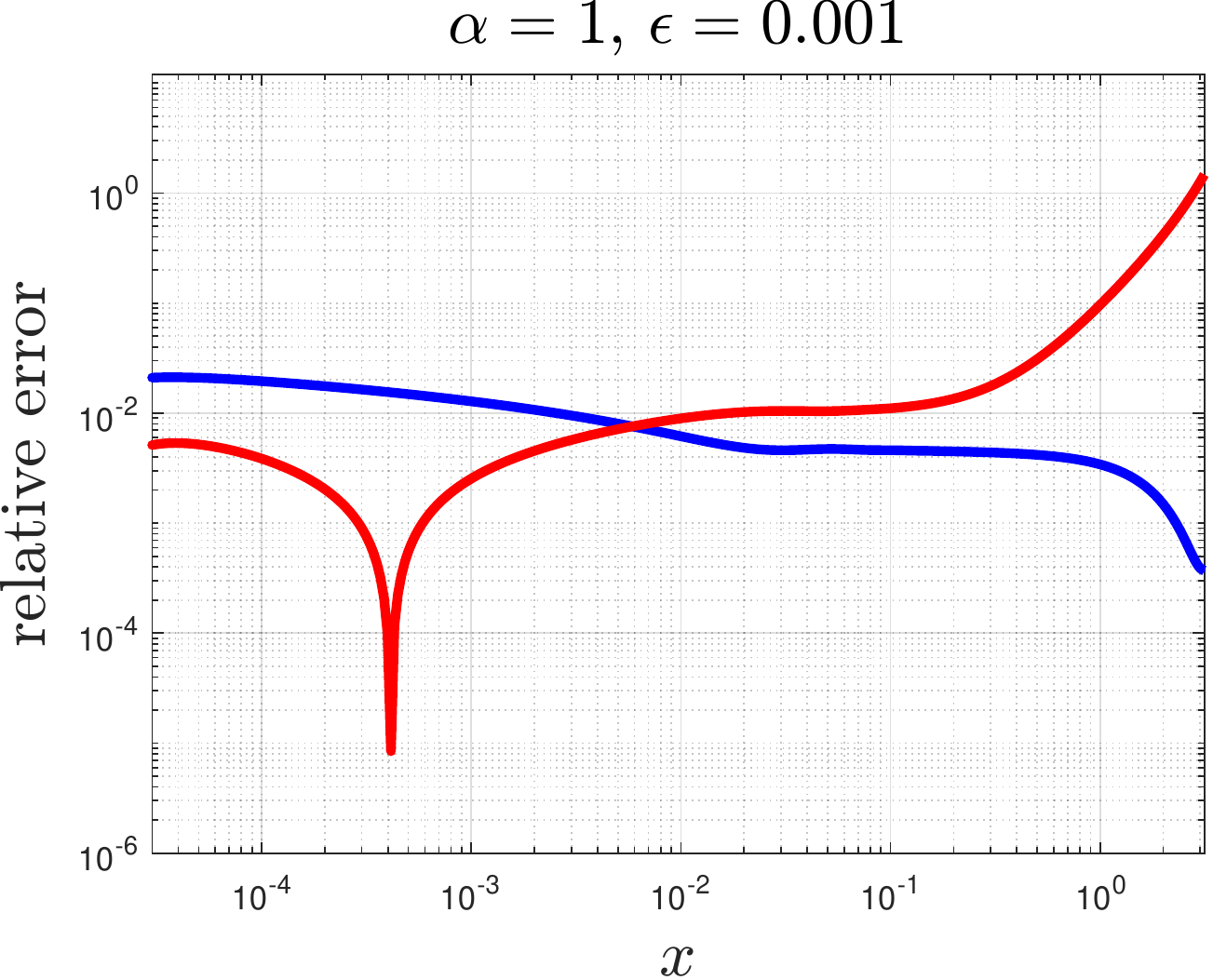}
	}
	\caption{The relative error of the asymptotic approximations (\ref{eq:blowupprofsimp}) (blue) and (\ref{eq:BUexpsmallx}) (red) at the blow-up time for small $x$. } \label{fig:smallxasympt}
\end{figure}

At the blow-up time, for $x$ exponentially small with respect to $\epsilon$, (\ref{eq:blowupprofsimp})  is no longer valid (see the discussion below (\ref{eq:wx0BU})). Instead,  
\begin{equation}
v \sim \frac{\epsilon\,e^{-\alpha}x^2}{2 - 8 \epsilon\, e^{-\alpha}\log(x^2)}, \label{eq:BUexpsmallx}
\end{equation}
which follows from (\ref{eq:wblowupprof}) and (\ref{eq:Fourbetagamma}). Figure~\ref{fig:smallxasympt} confirms that  (\ref{eq:BUexpsmallx}) has better accuracy than (\ref{eq:blowupprofsimp})  for small\footnote{Ideally, we would show the accuracy of (\ref{eq:BUexpsmallx}) for even smaller values of $x$ than those in Figure~\ref{fig:smallxasympt}. However, this would require high-precision computations. In standard double precision (with a machine precision of approximately $10^{-16}$) we cannot compute the errors for smaller $x$ because round-off errors prevent the accurate evaluation of $v$ from its numerically computed Fourier coefficients.} $x$ at the blow-up time.



 Regarding the strength of the singularity at $x = 0$, 
 it follows from (\ref{eq:BUexpsmallx}) that $v^{(n)} \sim 0$, $x \to 0$ for $n = 0, 1, 2$ but the third derivative, however, blows up as $v^{(3)} = O\left(\frac{1}{x\log^2 |x|}\right)$, $x \to 0$. A singularity with two bounded derivatives is consistent with the 
 $O(k^{-3})$ decay of the Fourier coefficients that was derived from (\ref{eq:blowupprofsimp})\footnote{The `global' blow-up profile (\ref{eq:blowupprofsimp})  has a stronger singularity at $x=0$ than the `local' blow-up profile (\ref{eq:BUexpsmallx}) since its second derivative blows up logarithmically at $x=0$ while for the local profile $v^{(2)} \sim 0$, $x \to 0$. The Fourier coefficients of the `global' profile nevertheless decay at the correct $O( k^{-3})$ rate since its second derivative, though unbounded, is integrable.}. In fact, the asymptotic approximation (\ref{eq:BUexpsmallx}) suggests that  the Fourier coefficients of the solution precisely at the blow-up time decay slightly faster than 
 the $O(k^{-3})$ suggested by~(\ref{eq:BUFourcoeffs}).  
 To  show this, symmetry and
 integration by parts can be used to obtain
 \[
 c_k(t_c)   = \frac{1}{\pi}\int_{0}^{\pi} v \cos k x = 
  \frac{1}{\pi k^3}\int_{0}^{\pi} v''' \sin k x\, dx.
 \]
 From (\ref{eq:BUexpsmallx}) it follows that
 $v \sim -x^2/(8\log( x^2 ))$, $v''' \sim 1/(8x\log^2(x))$,
 and by making the change of variable $x = \mu /k$ with $k \gg \mu$, 
 one obtains, 
 \begin{equation}
   c_k(t_c)  \sim   \frac{1}{8\pi k^3 \log^2 k} \int_{0}^{\infty} \frac{\sin \mu}{\mu} d\mu   =  \frac{1}{16 k^3 \log^2 k}. 
   \label{eq:FourBUc}
 \end{equation}
 The solid line in Figure~\ref{fig:smallxasympt} shows the estimate (\ref{eq:FourBUc}), which, unlike the estimate (\ref{eq:BUFourcoeffs}), is independent of $\epsilon$ and $\alpha$. The numerical Fourier coefficients shown in Figure~\ref{fig:smallxasympt} do not decay as fast as (\ref{eq:FourBUc}), which suggests that one would need to perform high-precision computations (to compute the solution closer to the blow-up time or with a larger number of Fourier coefficients) to observe the rate of decay predicted by (\ref{eq:FourBUc}).

For the possible continuation beyond blow-up one needs to confirm that the nonlinear term in the $v$-equation (\ref{eq:v}), i.e., $(v_x)^2/v$, remains bounded at $x = 0$ at the blow-up time. In fact, from (\ref{eq:BUexpsmallx}) it follows that 
the nonlinear term vanishes, according to  $ O(1/\log|x|)$, for $x \to 0$.

\section{Integrating through the singularity}\label{sec:postblowup}

The fact that the nonlinear term in~(\ref{eq:v}) remains bounded as $v$
approaches zero raises the intriguing possibility of numerically integrating
through the blow up.  This was tried in~\cite{Keller1993}, but the authors
found ``...the calculation actually continues the solution
slightly beyond the blow-up time of the original solution $u$. The method
becomes unstable a short time after the blow-up happens, however.''  
We report here on renewed efforts in this direction.
(Details of how our numerical strategies differ
from those of~\cite{Keller1993} are postponed to the end of the section.) 

\begin{figure}[htb]  
	\begin{center}
		\includegraphics[width=1\textwidth]{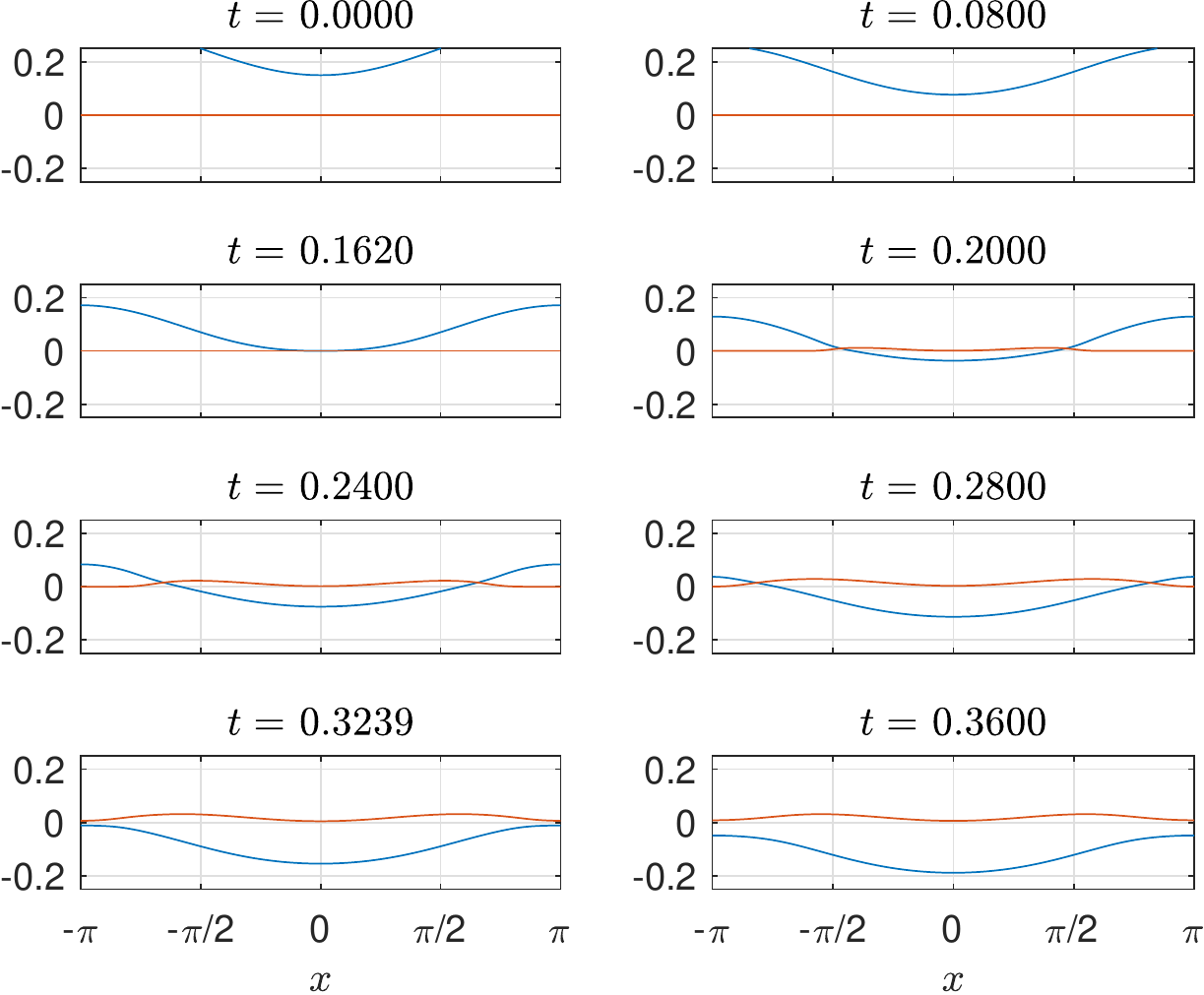}	
	\end{center}
	\caption{Solution to the $v$-equation~(\ref{eq:v}) corresponding to the initial 
		condition~(\ref{eq:ic}) with
		$\alpha = 0.25$, $\epsilon = 0.1$.  (The corresponding $u$ solution
		is shown in~Figure~\ref{fig:PostB}, and a movie of the dynamics
		can be seen in the supplementary material that accompanies
		this paper.)  The two noteworthy
		times are $t = t_c$ (third frame) and $t=2t_c$ (seventh frame), where 
	   $t_c=0.162$,  approximately.   At $t=t_c$ there is a zero of  $v$   at $x=0$,
	   after which the solution turns complex:   the real part is shown in blue and the
	   imaginary part in red. 
	   At $t=2t_c$ there are approximate zeros of  $v$ near  $x = \pm \pi$.  
	The solution shown is not unique (its conjugate is equally probable,
    as discussed
	in the text.) 
	}
	\label{fig:PostA}
\end{figure}

Figure~\ref{fig:PostA} shows a series of snapshots of the solution $v$ at different
times.   First, observe that there is no 
sign of instability as the solution passes smoothly through $v=0$ (third frame). 
What happens next might be unexpected, namely, the solution turns complex.
In addition, uniqueness is lost:
the complex conjugate of the solution shown is equally likely to appear
when parameters (such as the error tolerances in the time integration scheme or the number
of terms retained in the Fourier series~(\ref{eq:fourier1}))
are adjusted slightly.  Surprising as these results may be, both the non-uniqueness and the fact that the solution turns complex are consistent with theoretical results of~\cite{Masuda84}.

\begin{figure}[htb]  
	\begin{center}
  \includegraphics[width=1\textwidth]{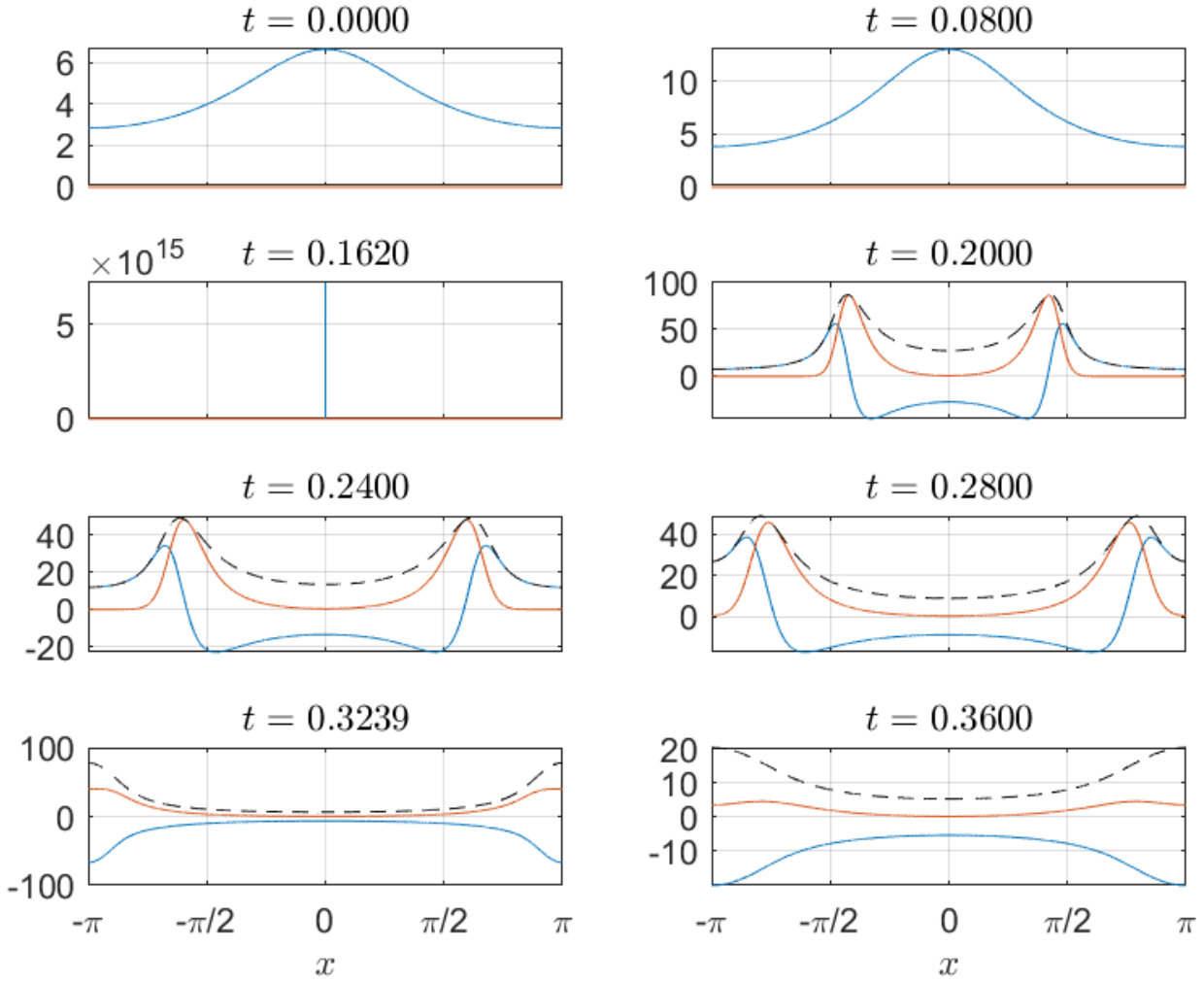}
	\end{center}
	\caption{Same as Figure~\ref{fig:PostA} but here the $u$-solution is shown.
		Note the blow up at $t=t_c$ and the near blow up at $t=2t_c$.
	The dashed curve is the modulus $|u|$ (shown only after the first blow up).  Note that
the scales on the vertical axes are different in each frame.
}
	\label{fig:PostB}
\end{figure}

Figure~\ref{fig:PostB} shows the $u$-solution that corresponds to the $v$-solution
in Figure~\ref{fig:PostA}.   After the solution turns complex at the
critical time (third frame)
the modulus of $u$ shows wave-like behaviour, with
two waves travelling in opposite directions from the origin until they
reach the edge of the domain $x = \pm \pi$.   Because of periodicity they meet up 
with similar waves 
from adjacent intervals and a second blow up almost
occurs, this time at $x = \pm \pi$  
(seventh frame).   The modulus grows considerably 
and we conjecture that by varying
the parameters in the initial condition 
a proper secondary blow up may be found. 
Computing over a longer time interval suggests that 
$u$ asymptotically 
approaches the constant (and real) solution $u \sim -1/t$ as $t \to \infty$.

We present the results of Figures~\ref{fig:PostA} and~\ref{fig:PostB} knowing full
well that we have little theory to draw on as validation.   The theoretical results of~\cite{Masuda84}
suggest the possibility of a complex and non-unique solution post blow-up, but does not
allow for a quantitative comparison.   
Nevertheless, the following heuristic 
observations support the validity of the results shown here. 

\begin{figure}[h!]  
	\begin{center}
		\includegraphics[width=1\textwidth]{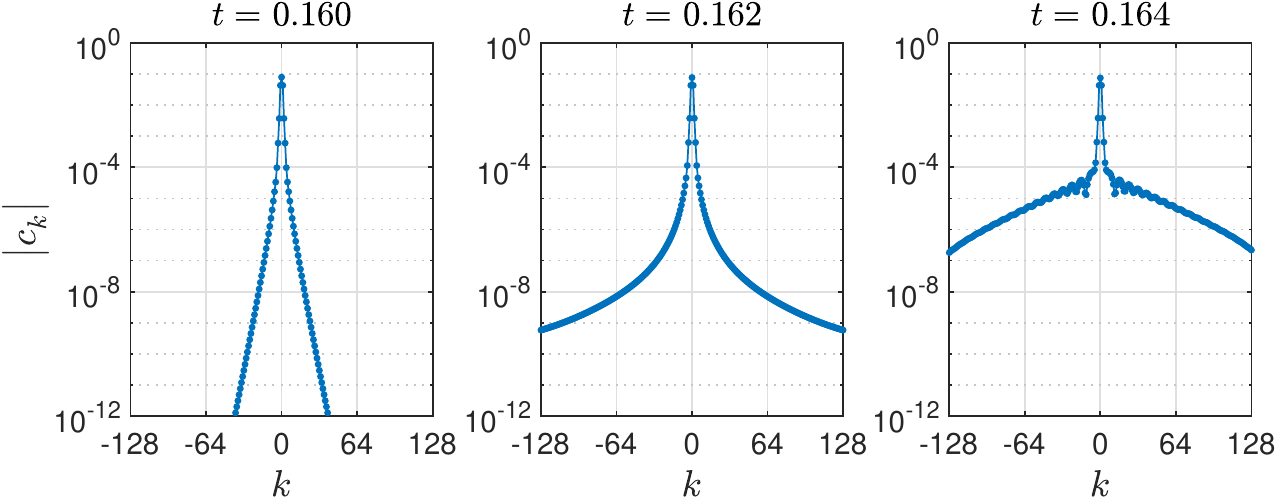}	
	\end{center}
	\caption{Modulus of the Fourier coefficients $c_k$
		of the solution $v$ shown in Figure~\ref{fig:PostA}. 
		The times are  just prior to $t=t_c$, at  $t=t_c$, and
		shortly after.  (The value of $t_c$   
		is given to more digits in Table~\ref{tab:tc}.)
		In the first frame the coefficients decay exponentially,
		indicative of a function analytic in a neighbourhood of the real axis.   At $t=t_c$ the
		coefficients decay algebraically 
		as given by~(\ref{eq:FourBUc}), 
		because of the influence of the singularity of $v$
		as it reaches the origin.  Immediately afterwards  
		the exponential decay is recovered. 
}  
	\label{fig:Fourier}
\end{figure}

Firstly, the fact that the numerical solution approaches
$-1/t$ as $t \to \infty$ provides some confidence, as this 
solves~(\ref{eq:u}) exactly.    Secondly, as noted before, the singularity in the $v$-equation at the critical time is rather weak,
leading to Fourier coefficients that decay 
at the relatively rapid rate 
$c_k = O(1/(k^3 \log^2 k))$; recall~(\ref{eq:FourBUc}). 
While still much slower than
the typical exponential decay rate for analytic periodic functions, we conjecture that this decay is 
nevertheless sufficiently rapid
that the accuracy loss of the spectral method at the critical time is not disastrous. 
Figure~\ref{fig:Fourier} shows the Fourier coefficients near and at the critical time. 

As far as we know, the only published results that
deal with numerical computations of  post-blow-up solutions are~\cite{Cho16,Takayasu22}. 
These authors based their computations on
complexification of the $t$-variable.
Following this idea we integrated along a path that
contains a semi-circle in the complex $t$-plane, 
centred at the estimated singularity.   Using this 
approach we obtained the same results as those shown
in Figures~\ref{fig:PostA} and~\ref{fig:PostB}. 

It should be noted that the methods of~\cite{Cho16,Takayasu22} are based on the $u$-equation, not the $v$-equation.   Because the solution to the $u$ equation grows without bound 
near the critical time, these methods cannot compute  blow-up solutions 
in the immediate neighbourhood of the critical time accurately.  For the same reason, 
estimates of blow-up times  
are unreliable.   Adapting the methods of~\cite{Cho16,Takayasu22} from
the $u$-equation to the $v$-equation might be a
worthwhile future project. 


The observation of a complex solution after the critical time has lead
us to conjecture that the failure of the 
finite difference method of~\cite{Keller1993} was because it was
almost surely  coded as a real system using real arithmetic.  
One way to allow for possible complex solutions is to use 
complex arithmetic (default in \uppercase{MATLAB}) or simply
to split the $v$-equation~(\ref{eq:v}) into its
real and imaginary parts.   We integrated such a separated system with
the Fourier spectral method mentioned in the first section,
which is how the results of 
Figures~\ref{fig:PostA}--\ref{fig:Fourier} were computed.
The same idea applies, however, to the finite difference method of~\cite{Keller1993}.

Initialising the imaginary part to strictly zero values, 
the ODE software signals a singularity at the critical time.  
Initializing it with a small, random perturbation on the order of 
machine roundoff level ($\sim 10^{-16}$), however, allowed the solution
$v$ to pass through the critical $v=0$  without any sign of instability, as shown in 
Figure~\ref{fig:PostA}.


One priority for future investigations is a better understanding of the transition from a real to a complex solution and the associated
non-uniqueness.   In our case the nonzero imaginary part is
triggered by noise at the level of roundoff error.  The radomness
dictates whether the continuation
is with one solution or with its complex conjugate.








\section{Conclusions}  

We investigated, asymptotically and numerically, point blow-up solutions to a periodic nonlinear heat equation (\ref{eq:u})  with nearly flat initial data by considering the solution $u$ in the reciprocal variable $v = 1/u$.  We derived asymptotic approximations for the solution on the entire spatial interval and from $t=0$ up to and including at  the blow-up time, for which we also derived a second-order approximation. Due to the high accuracy of the Fourier spectral method (including at the blow-up time at which the $v$-solution has a weak singularity, unlike the $u$-solution), we were able to check numerically the validity of the asymptotics. We believe it is unusual for a single numerical method to confirm asymptotics in multiple regimes since typically numerical methods are weak, or highly inefficient, in most asymptotic limits unless they are highly specialised.   
The key to the success of the numerical method used here is the
fact that it approximates the relatively well-behaved 
$v$-equation~(\ref{eq:v}) rather than the $u$-equation~(\ref{eq:u})
whose solution becomes unbounded. 

The investigations in this paper point to a number of topics for future research, not least of which is the validity of the post-blow-up solutions computed in section~\ref{sec:postblowup}. 
In addition, the dynamics of  complex singularities of blow-up solutions to (\ref{eq:u}) for a larger class of initial conditions, including initial data leading to non-generic forms of blow up, will be investigated  in~\cite{FKW}, also via a combination of asymptotic and numerical methods. In~\cite{FKW} we shall also explore  the singularity structure of these blow-up solutions on their Riemann surfaces in the complex $x$-plane.







\begin{acknowledgements}
The first author is grateful to Saleh Tanveer for stimulating and insightful discussions.
	This research was started while the three authors were in residence
	 at the Isaac Newton Institute for Mathematical Sciences as part of the  programme 
	{\sl Complex analysis: techniques, applications and computations}.
	This programme was supported by: EPSRC grant number EP/R014604/1. The work of the first author was also supported by the Leverhulme Trust Research Project Grant RPG-2019-144.
	A grant to the third author from the H.B.~Thom foundation of Stellenbosch University is also
	gratefully acknowledged.    
\end{acknowledgements}  
 
\bibliographystyle{plain}
\bibliography{BlowUp}

\appendix
 \section{Analysis via the method of matched asymptotic expansions}\label{sec:matchedasympt}
 
 
 


\subsection{Truncated Fourier expansion}\label{sect:tFourexp}

In this section we revisit and expand the
analysis based on the two-mode Fourier truncation outlined in section~\ref{sec:twomode}.
Recall that by substituting~(\ref{eq:ansatz}) into the $v$-equation (\ref{eq:v}) 
and neglecting the $\cos 2x$ term the system~(\ref{eq:twomode}) was obtained.  From this, one finds by a self consistency argument that near blow-up, generically  (i.e.~even for $\epsilon = O(1)$)

\begin{equation}
a \sim a_c + (1 + 2a_c)(t_c-t), \qquad b \sim a_c - a_c(t_c-t), \qquad t \to t_c^{-},  \label{eq:FourBUgenab}
\end{equation}
where $a_c$ and the blow up time $t_c$ are positive constants, so that the blow-up behaviour associated with (\ref{eq:ansatz}) takes the form
\begin{equation}
v \sim (1 + 3a_c)(t_c-t) + \frac{1}{2}a_c x^2, \qquad t \to t_c^{-}, \qquad x = O\left((t_c-t)^{1/2}\right).  \label{eq:FourBUgenv}
\end{equation}
We emphasise that (\ref{eq:FourBUgenab})--(\ref{eq:FourBUgenv}) describe the blow-up behaviour of  (\ref{eq:twomode}) rather than that of (\ref{eq:v}): (\ref{eq:FourBUgenv}) has features in common with, but is not a valid representation of, the blow-up behaviour of the full PDE (\ref{eq:v}), a key point to which we shall return.

In keeping with our goal of characterising the behaviour of (\ref{eq:v}) for near-flat initial conditions, we now develop a fully analytic asymptotic description of the behaviour of the two-mode system (\ref{eq:twomode}) for initial conditions (\ref{eq:FICs}) with $\alpha = O(1)$. There are two timescales; the first coincides with (\ref{eq:F1sttimescale}), so that at leading order (\ref{eq:a0s}) and (\ref{eq:b0s}) follow.

On the second time scale, and near blow up, we set
\begin{equation}
t = \alpha + \epsilon T,\quad b = \epsilon B, \quad a = \epsilon A_0(T) + O(\epsilon^2),   \quad 
B =  B_0(T) + O(\epsilon),  \label{eq:Ft2lims}
\end{equation}
 with $T = O(1)$, so that
\begin{equation}
\begin{split}
A_0\frac{dA_0}{dT} + \frac{1}{2}B_0  \frac{dB_0}{dT} & = -A_0, \\
B_0\frac{dA_0}{dT} + A_0  \frac{dB_0}{dT} & = -B_0.
\end{split} \label{eq:A0B0}
\end{equation}
 Since $A_0 \pm \frac{1}{\sqrt{2}}B_0$ each satisfy
\begin{equation*}
\Phi\frac{d\Phi}{dT} = - \Phi
\end{equation*}
we find on matching with (\ref{eq:a0s})--(\ref{eq:b0s})  that
\begin{equation}
A_0 = - T, \qquad B_0 = e^{-\alpha}  \label{eq:A0B0s}
\end{equation}
(a result that relies on the observation that $a$ has no $O(\epsilon)$ contribution for $t = O(1)$, see (\ref{eq:b0s})). Thus
\begin{equation*}
t_c \sim \alpha - \epsilon e^{-\alpha}, \qquad a_c \sim \epsilon e^{-\alpha}, \qquad \epsilon \to 0^{+}
\end{equation*}
in (\ref{eq:FourBUgenab})--(\ref{eq:FourBUgenv}). Moreover, while the rescaling (\ref{eq:Ft2lims}) leads to a modified balance, i.e. (\ref{eq:A0B0}) in place of (\ref{eq:a0s}) and (\ref{eq:b0s}), the result (\ref{eq:A0B0s}) implies that the solution passes unscathed through $T = O(1)$, with (\ref{eq:a0s})--(\ref{eq:b0s}) being recovered for $t > \alpha$. Thus with the two-mode approximation (\ref{eq:ansatz}) and (\ref{eq:twomode1}) continuation through blow up seems to be straightforward, in contrast to that of the $v$-equation (\ref{eq:v}), as we shall subsequently demonstrate.

\subsection{Truncated Taylor expansion}\label{sect:Texp}

Much of what follows in this subsection revisits results from~\cite{Keller1993}, which we derive using a different  approach (namely, matched asymptotic expansions).

As in section~\ref{sec:twomode}, we also consider the truncated Taylor approximation (\ref{eq:Taylans}), which gives rise to the system (\ref{eq:twomodeKL}) whose  blow-up behaviour takes the form (again, by a self-consistency argument)
\begin{equation}
a \sim t_c - t, \qquad b \sim \frac{1}{8(-\log(t_c-t) + b_c)}, \qquad t \to t_c^{-}  \label{eq:TgenBU}
\end{equation}
for constants $t_c$ and $b_c$. In contrast to the results of the previous subsection, (\ref{eq:twomodeKL}) \emph{does} capture the blow-up behaviour of the full PDE (\ref{eq:v}), a point to which we shall also come back.

We now return to initial conditions of the form (\ref{eq:FICs}), with the scalings (\ref{eq:F1sttimescale}) applying for $t = O(1)$, so expanding in the form
\begin{equation*}
a \sim a_0 + \epsilon a_1, \qquad b = \epsilon B, \qquad B \sim B_0 + \epsilon B_1,
\end{equation*}
(\ref{eq:twomodeKL})  implies
\begin{equation}
a_0 = \alpha - t, \quad B_0 = 1,\quad a_1 = 2t, \quad B_1 = 8\log\left( (\alpha - t)/\alpha  \right).  \label{eq:Ta0b0s}
\end{equation}

Under the scalings (\ref{eq:Ft2lims}) the leading-order balances do not change and the results of relevance below read
\begin{equation}
A_0 = 2\alpha - T, \quad B_0 = 1, \quad B_1 = 8\log\left((2\alpha - T)/\alpha  \right),  \label{eq:Tt2sols}
\end{equation}
the matching into $a_1$ in (\ref{eq:Ta0b0s}) leading to the $2\alpha$ contributions. The final scale in this case is then more subtle than those above: we set
\begin{equation}
a = (T_c - T)\hat{a}(\tau), \quad B = \hat{b}(\tau), \quad \tau = - \epsilon\log(T_c - T),  \label{eq:Texpt}
\end{equation}
where 
\begin{equation}
T_c \sim 2\alpha, \qquad t_c \sim \alpha + 2\epsilon\alpha, \qquad \epsilon \to 0^{+}. \label{eq:TaylBUt}
\end{equation}
The introduction of $\tau$ in (\ref{eq:Texpt}) is associated with the expansion for $B$ in (\ref{eq:Tt2sols}) disordering and (\ref{eq:twomodeKL}) becomes
\begin{equation}
\begin{split}
& \epsilon\frac{d\hat{a}}{d\tau} - \hat{a} = 2\epsilon \hat{b} - 1,  \\
& \hat{a} \frac{d\hat{b}}{d\tau} = -8\hat{b}^2,
\end{split}  \label{eq:Ttwosysh}
\end{equation}
so that, on matching with (\ref{eq:Tt2sols}),
\begin{equation}
\hat{a}_0 = 1, \qquad \hat{b}_0 = \frac{1}{1 + 8\tau}, \qquad \hat{a}_1 = -\frac{2}{1+8\tau},  \label{eq:t3sols}
\end{equation}
so that $b_c \sim 1/(8\epsilon)$ in (\ref{eq:TgenBU}).

For reasons that will become apparent below, it is instructive to record the implications of (\ref{eq:t3sols}) under the scaling
\begin{equation}
x = \epsilon^{1/2}X  \label{eq:xX}
\end{equation}
whereby
\begin{equation*}
v \sim \epsilon\left( T_c - T + \frac{\epsilon}{1 + 8\tau}\left(X^2 - 2 (T_c - T)  \right)  \right);
\end{equation*}
this will reappear in the analysis below in describing the local blow-up behaviour, but also confirms that (\ref{eq:Taylans}) cannot describe the spatial profile at blow up (i.e. at $\tau = +\infty$).

\subsection{More general near-flat initial data} \label{sec:appgenflat}

We now turn to the derivation of the blow-up behaviour, subjecting (\ref{eq:v}), the full PDE, to the initial data
\begin{equation*}
v(x,0) = \alpha + \epsilon V(x)
\end{equation*}
with, again, $0 < \epsilon \ll 1$. It is striking that this limit allows a near complete analytical description of the transition to blow up through the fully nonlinear regime. Three time scales are required. On the first time scale:
\begin{equation}
t = O(1), \qquad v \sim v_0(t) + \epsilon v_1(x,t) + \epsilon^2 v_2(x,t).  \tag{I}
\end{equation}
All three terms in this expansion are required for what follows. An immediate result is that
\begin{equation*}
v_0 = \alpha - t
\end{equation*}
leading to 
\begin{equation}
\frac{\partial v_1}{\partial t} = \frac{\partial^2 v_1}{\partial x^2}, \qquad v_1(x,0) = V(x).  \label{eq:heatv1}
\end{equation}
We shall consider general $V(x)$, but in the case of real-valued 
$2\pi$-periodic initial conditions (or zero Neumann boundary conditions on a finite domain), we can Fourier decompose in the usual way, so that
\begin{equation*}
v_1(x,t) =  \sum_{k=-\infty}^{\infty}\ a_k e^{ikx-k^2 t}.
\end{equation*}
We note this special case for two reasons -- firstly for its relevance to section~\ref{sect:tFourexp} and secondly for the obvious observation that the high-frequency modes are rapidly decaying, providing additional motivation for the analysis of sections~\ref{sec:twomode} and~\ref{sect:tFourexp} but being also in some respects counter-intuitive, given that point blow up subsequently ensues. 

We define 
\begin{equation*}
\Phi_1(x) = v_1(x,\alpha)
\end{equation*}
so that
\begin{equation}
v_1(x,t) \sim \Phi_1(x) - (\alpha - t)\Phi_1''(x), \qquad t \to \alpha^{-}.  \label{eq:phi1def}
\end{equation}
At next order we have
\begin{equation}
\frac{\partial v_2}{\partial t} - \frac{\partial^2 v_2}{\partial x^2} = -\frac{2}{\alpha - t}\left(\frac{\partial v_1}{\partial x}\right)^2, \qquad v_2(x,0) = 0,  \label{eq:v2eq}
\end{equation}
so that
\begin{equation}
\begin{split}
v_2(x,t) \sim & 2\left(\Phi_1'(x)\right)^2\log(\alpha - t) + \Phi_2(x) - 2\left[\left( \Phi_1' \right)^2  \right]''(\alpha - t)\log(\alpha - t)  \\
& + \left( 4\left(\Phi_1''  \right)^2 - \Phi_2''  \right)(\alpha - t), \qquad t \to \alpha^{-},
\end{split}  \label{eq:phi2def}
\end{equation}
where (\ref{eq:phi2def}) serves to define $\Phi_2(x)$.

On the second time scale:
\begin{equation}
v = \epsilon w, \qquad w \sim w_0(x,T) + \epsilon\log(1/\epsilon) w_1(x,T) + \epsilon w_2(x,T).  \tag{II}
\end{equation}
Here $t = t_c(\epsilon) + \epsilon T$, $T =  O(1)$ with $t_c(0) = \alpha$ and with correction terms to $t_c$ identified below\footnote{The notation here differs somewhat from that of the previous subsections.}. We shall also need to consider the rescaling (\ref{eq:xX}). We have
\begin{equation*}
w \frac{\partial w}{\partial T} = \epsilon\left(w \frac{\partial^2 w}{\partial x^2} - 2 \left( \frac{\partial w}{\partial x} \right)^2 \right) - w
\end{equation*}
so, matching with (\ref{eq:phi1def}), (\ref{eq:phi2def}), and defining $\beta_1$ and $\beta_2$ via
\begin{equation}
t_c(\epsilon) \sim \alpha + \epsilon \beta_1 + \epsilon^2 \beta_2, \qquad \epsilon \to 0 \label{eq:beta12}
\end{equation}
it follows for $x = O(1)$ that
\begin{equation}
w_0 = -T + \Phi_1(x) - \beta_1, \label{eq:w0s}
\end{equation}
\begin{equation}
w_1 = -2\left(  \Phi_1'(x) \right)^2 \label{eq:w1s}
\end{equation}
and hence
\begin{equation*}
\frac{\partial w_2}{\partial T} = \Phi_1''(x) - \frac{2}{-T+\Phi_1(x)-\beta_1}\left(\Phi_1'(x)\right)^2,
\end{equation*}
implying
\begin{equation}
w_2 = - (-T)\Phi_1''(x) + 2\left(\Phi_1'(x)\right)^2\log\left( -T+\Phi_1(x) - \beta_1 \right) + \Phi_2(x) - \beta_2.\label{eq:w2s}
\end{equation}
We take blow up to occur at $x=0$, requiring that $\Phi_1'(0) = \Phi_2'(0) = 0$ and take
\begin{equation}
\Phi_1(x) \sim \beta_1 + \gamma_1 x^2, \qquad
\Phi_2(x) \sim \beta_2 + \gamma_2 x^2, \qquad
x \to 0  \label{eq:phi12}
\end{equation}
for constants $\beta_{1,2}$, $\gamma_{1,2}$ with $\gamma_1 >0$\footnote{The case $\gamma_1=0$ corresponds to non-generic forms of blow up -- we shall not pursue such matters here.} and where the requirement that $w=0$ at $x=0$, $T=0$ implies that $\beta_{1,2}$ in (\ref{eq:beta12}) are specified by (\ref{eq:phi12}), with $\Phi_{1,2}$ being determined by the linear problems (\ref{eq:heatv1}) and (\ref{eq:v2eq}). 

For $X = O(1)$, we first generate the required matching conditions from (\ref{eq:phi1def}), (\ref{eq:phi2def}) and (\ref{eq:phi12}), whereby
\begin{equation}\label{eq:wexp}
\begin{split}
w \sim & -T + \left(\epsilon\gamma_1 - 8\epsilon^2\log(1/\epsilon)\gamma_1^2 + 8\epsilon^2\gamma_1^2\log(-T) + \epsilon^2\gamma_2   \right)\left(X^2 - 2(-T)  \right) \\
& + 16\epsilon^2\gamma_1^2(-T),
\end{split}
\end{equation}
wherein we have retained only the required terms in $X^0$ and $X^2$ -- in general $v_1$ will also lead to contributions of the form
\begin{equation*}
\epsilon^{3/2}\left( X^3 - 6(-T)X  \right), \qquad
\epsilon^2\left( X^4 - 12(-T)X^2 + 12(-T)^2   \right)
\end{equation*}
in $w$ but these can be neglected for our purposes, being sub-dominant as $T \to 0^{-}$ with $X = O\left( (-T)^{1/2} \right)$.

We have
\begin{equation}
w\frac{\partial w}{\partial T} = w\frac{\partial^2 w}{\partial X^2} - 2\left(\frac{\partial w}{\partial X}\right)^2 - w  \label{eq:wXeq}
\end{equation} 
and from this we find that the relevant terms for $T = O(1)$ simply reproduce the matching condition (\ref{eq:wexp}) obtained from the expansion for $t = O(1)$. Importantly, the expansion (\ref{eq:wexp}) disorders for $\epsilon\log(-T) = O(1)$, as in section~\ref{sect:Texp}, leading us on to our final time scale, as follows:
\begin{equation}
  \tau = O(1), \quad w = (-T)g, \quad g \sim 1 + \epsilon g_1(\xi, \tau) + \epsilon^2\log(1/\epsilon)g_2(\xi,\tau) + \epsilon^2 g_3(\xi, \tau)  \tag{III}
\end{equation}
where
\begin{equation}
\tau = -\epsilon\log(-T), \qquad \xi = \frac{X}{(-T)^{1/2}}.
\label{eq:solanz}
\end{equation}

We emphasise that much of what follows reconstructs known blow-up behaviour (see
the review of~\cite{Galaktionov2002}, for example.) The novelty here lies in the focus
on the near-flat initial data, which allows more detailed characterization of the
full spatial behaviour.

The near self-similar solution ansatz~(\ref{eq:solanz}) transforms (\ref{eq:wXeq}) to
\begin{equation*}
\epsilon g \frac{\partial g}{\partial \tau} - g^2 + \frac{1}{2}\xi g\frac{\partial g}{\partial \xi} = g \frac{\partial^2 g}{\partial \xi^2} - 2\left( \frac{\partial g}{\partial \xi} \right)^2 - g
\end{equation*}
so that
\begin{equation}
-g_1 + \frac{1}{2}\xi \frac{\partial g_1}{\partial \xi} = \frac{\partial^2 g_1}{\partial \xi^2},\label{eq:g1eq}
\end{equation}
implying that
\begin{equation}
g_1 = \sigma_1(\tau)\left(\xi^2 - 2\right);\label{eq:g1sol}
\end{equation}
similarly
\begin{equation}
g_2 = \sigma_2(\tau)\left(\xi^2 - 2\right) \label{eq:g2sol}
\end{equation}
but $g_3$ satisfies (following cancellation of a number of terms using (\ref{eq:g1eq}))
\begin{equation}
\frac{d\sigma_1}{d\tau}\left(\xi^2 - 2\right) - g_3 + \frac{1}{2}\xi \frac{dg_3}{d\xi} = \frac{d^2g_3}{d\xi^2} - 8\sigma_1^2\xi^2.  \label{eq:g3eq}
\end{equation}
Since we need to preclude exponential growth of $g_3$ (i.e. to exclude contributions with $\log g_3 \sim \xi^2/4$ as $\xi \to \pm \infty$), (\ref{eq:g3eq}) both requires that 
\begin{equation}
g_3 = -2\frac{d\sigma_1}{d\tau} + \sigma_3(\tau)\left( \xi^2 - 2\right)  \label{eq:g3sol}
\end{equation}
and generates the solvability condition
\begin{equation}
\frac{d\sigma_1}{d\tau} = -8\sigma_1^2;\label{eq:sigma1solv}
\end{equation}
that the $O(\epsilon)$ term in $g$ is only fully determined via (\ref{eq:g3eq}) necessitates that the expansion be taken up to $O(\epsilon^2)$. Matching with (\ref{eq:wexp}) then requires that
\begin{equation}
\sigma_1 = \frac{\gamma_1}{1 + 8\gamma_1 \tau}, \qquad
-2\frac{d\sigma_1}{d\tau} = \frac{16\gamma_1^2}{(1 + 8\gamma_1\tau)^2}, \label{eq:sig1sol}
\end{equation}
simultaneously confirming matching with the second, fourth and final terms in (\ref{eq:wexp}). The calculation of $\sigma_2$ and $\sigma_3$ requires solvability conditions at yet higher orders (which we will not pursue), with (\ref{eq:wexp}) requiring that
\begin{equation*}
\sigma_2(0) = -8\gamma_1^2, \qquad \sigma_3(0) = \gamma_2.
\end{equation*}

By analogy with the subdivision into $x$ and $X$ in (II), we need also to consider $\eta = O(1)$, where
\begin{equation*}
\eta = \frac{x}{(-T)^{1/2}}, \qquad \eta = \epsilon^{1/2}\xi,
\end{equation*}
though a third scale (namely $x = O(1)$) will also require consideration here. Since
\begin{equation*}
\epsilon g \frac{\partial g}{\partial \tau} - g^2 + \frac{1}{2}\eta g \frac{\partial g}{\partial \eta} = \epsilon\left(g \frac{\partial^2 g}{\partial \eta^2} - 2\left(\frac{\partial g}{\partial \eta}  \right)^2  \right) - g,
\end{equation*}
setting
\begin{equation*}
g \sim G_0(\eta, \tau) + \epsilon\log(1/\epsilon)G_1(\eta,\tau) + \epsilon G_2(\eta,\tau)
\end{equation*}
gives on matching into $\xi = O(1)$
\begin{equation}
G_0 = 1 + \sigma_1(\tau)\eta^2, \qquad G_1 = \sigma_2(\tau)\eta^2  \label{eq:G0G1s}
\end{equation}
and
\begin{equation*}
\frac{1}{2}\eta \frac{dG_2}{d \eta} - G_2 = 2\sigma_1 - \frac{d\sigma_1}{d\tau}\eta^2 - \frac{8\sigma_1^2\eta^2}{1 + \sigma_1\eta^2},
\end{equation*}
so that, again matching into $\xi = O(1)$,
\begin{equation}
G_2 = -2\sigma_1(\tau) + 8\sigma_1^2(\tau)\eta^2\log\left(1 + \sigma_1(\tau)\eta^2\right) + \sigma_3(\tau)\eta^2;\label{eq:G2s}
\end{equation}
were $\sigma_1$ not a solution to (\ref{eq:sigma1solv}), a term in $\log \eta$ would also be present here, further clarifying the status of (\ref{eq:sigma1solv}) as a solvability condition. 

Finally, we can simply set $T=0$ in (\ref{eq:w0s})--(\ref{eq:w2s}), noting that $\tau = O(1)$ corresponds to exponentially small $T$, to obtain the profile at blow up for $x = O(1)$. Thus as $t \to t_c^{-}$, $\epsilon \to 0$ with $x = O(1)$ we have
\begin{equation}
\begin{split}
w \sim & \: \Phi_1(x) - \beta_1 - 2\epsilon\log(1/\epsilon)\left(\Phi_1'(x)\right)^2 \\
& + \epsilon\left(2\left(\Phi_1'(x)\right)^2\log(\Phi_1(x) - \beta_1) + \Phi_2(x) - \beta_2   \right).
\end{split} \label{eq:wBU}
\end{equation}
The expression (\ref{eq:wBU}) has small-$x$ behaviour
\begin{equation}
w  \sim \gamma_1 x^2 - 8\epsilon\log(1/\epsilon)\gamma_1^2x^2 + \epsilon\left( 8\gamma_1^2x^2\log(\gamma_1x^2) + \gamma_2x^2 \right), \label{eq:wx0BU}
\end{equation}
The expressions (\ref{eq:wBU})--(\ref{eq:wx0BU}) do not apply for exponentially small $x$, however. Instead, we need to extract from (\ref{eq:G0G1s})--(\ref{eq:G2s}) the terms relevant for large $\eta$, namely
\begin{equation*}
    g \sim \sigma_1 \eta^2 + 16 \epsilon \sigma_1^2 \eta^2 \log \eta,
\end{equation*}
these being of the same order for $\epsilon\log \eta = O(1)$.

Reconstructing $w$ from these using (\ref{eq:sig1sol}), we have
\begin{equation*}
    w \sim \frac{\gamma_1}{1 + 8\gamma_1 \tau}x^2 + \frac{16 \epsilon \gamma_1^2 x^2}{(1+8\gamma_1 \tau)^2}\log\eta,
\end{equation*}
which, using
\begin{equation*}
\tau = -2\epsilon\log x + 2\epsilon\log \eta,
\end{equation*} 
implies that
\begin{equation}
w \sim \frac{\gamma_1 x^2}{1 + 16\epsilon \gamma_1 \log(1/x)} \label{eq:wblowupprof}   
\end{equation}
describes the profile at blow up for $x$ exponentially small with respect to $\epsilon$ (having taken the various limits in the appropriate order); (\ref{eq:wblowupprof}) matches with the relevant terms in (\ref{eq:wx0BU}) for larger $x$.

We note that non-analytic (i.e. logarithmic) terms here occur as a matter of course, in $(t_c-t)$ in (\ref{eq:sig1sol}) and in $x$ in (\ref{eq:wx0BU}): these arise constructively in the current analysis rather than being introduced a priori as part of a solution ansatz. That the current limit provides a detailed asymptotic characterisation of the profile (\ref{eq:wBU}) at blow up for almost all $x$, not just close to blow-up point, is also noteworthy. 

\subsubsection{Comparison with truncated expansions}\label{sec:apptrunccomp}

The above systematic asymptotic analysis clarifies the extent of applicability of the ad-hoc approximations treated in sections~\ref{sect:tFourexp} and~\ref{sect:Texp}.

Starting with a comparison with the results of section~\ref{sect:Texp}, we have already noted that the truncated Taylor expansion correctly captures the blow-up behaviour, as can be substantiated by the following observations. For $V(x) = x^2$ it follows that 
\begin{equation}
\begin{split}
& v_1 = x^2 + 2t, \qquad v_2 = 8x^2\log\left( (\alpha - t)/\alpha  \right) - 16(\alpha - t)\log\left((\alpha - t)/\alpha\right) - 16t,   \\
& \Phi_1 = x^2 + 2\alpha, \qquad  \Phi_2 = -8x^2\log\alpha - 16\alpha,   
\end{split}\label{eq:phi12Tayl}
\end{equation}
which imply (see (\ref{eq:phi12}) and (\ref{eq:beta12}))
\begin{equation}
    \beta_1 = 2\alpha, \qquad \gamma_1 = 1, \qquad \beta_2 = -16\alpha, \qquad \gamma_2 = -8\log\alpha, \qquad t_c \sim \alpha + 2\alpha \epsilon - 16\alpha \epsilon^2.  \label{eq:Taylbetagamma}
\end{equation}
More importantly, (\ref{eq:wexp}), (\ref{eq:g1sol}), (\ref{eq:g2sol}), (\ref{eq:g3sol}) and (\ref{eq:G0G1s}) are all quadratic in $x$ for general $V(x)$, so the 
approximation in section~\ref{sect:Texp} represents an attractor in that sense. As already noted, it cannot capture the spatial blow-up profile, the first manifestation of this being represented by the logarithmic term in (\ref{eq:G2s}); see also  (\ref{eq:wx0BU}) and (\ref{eq:wblowupprof}).  It is especially noteworthy that (\ref{eq:g3sol}) is a quadratic in $\xi$; if the corresponding analysis is undertaken on (\ref{eq:u}) rather than on (\ref{eq:v}), a $\xi^4$ term arises at that order. This represents a hidden benefit of the $v$ formulation.

Turning now to section~\ref{sect:tFourexp}, corresponding to $V(x) = -\cos x$, we can both exemplify the form of the profile (\ref{eq:wBU}) at blow up and indicate where the analysis of section~\ref{sect:tFourexp} breaks down in describing the behaviour of the full PDE. In this case
\begin{equation}
\begin{split}
& v_1 = -e^{-t}\cos x, \quad v_2 = -\int_0^t\frac{e^{-2s}}{\alpha - s}ds - e^{-4t}\int_0^t\frac{e^{2s}}{\alpha - s}ds\cos(2x),  \ \\
& \Phi_1 = -e^{-\alpha}\cos x, \quad
 \Phi_2 = -C_1 - C_2  -\left(C_1 + C_3  \right)\cos 2x,
\end{split}\label{eq:phi12Four}
\end{equation}
where
\begin{equation}
C_1 = e^{-2\alpha}\log\alpha, \qquad 
C_2 = \int_0^{\alpha}\frac{e^{-2t}-e^{-2\alpha}}{\alpha - t}dt, \qquad
C_3 = e^{-4\alpha}\int_0^{\alpha}\frac{e^{2t}-e^{2\alpha}}{\alpha - t}dt,  \label{eq:cdefs}
\end{equation}
which imply (see (\ref{eq:phi12}))
\begin{equation}
    \beta_1 = -e^{-\alpha}, \qquad \gamma_1 = \frac{e^{-\alpha}}{2}, \qquad
    \beta_2 = -2C_1-C_2-C_3, \qquad
    \gamma_2 = 2(C_1+C_3).  \label{eq:Fourbetagamma}
\end{equation}
That $v_2$ and $\Phi_2$ contain $\cos 2x$ contributions is already indicative of the failure of the truncated Fourier expansion close to the blow up; nevertheless, $v_1$ and $w_0$ (in (\ref{eq:w0s})) provide the dominant spatial dependencies on the relevant scales and retain the form of section~\ref{sect:tFourexp} (and the $O(\epsilon)$ term in (\ref{eq:wexp})  can be viewed as being associated with the latter's Taylor expansion). The truncation 
is, however, entirely unable to reproduce the behavior for $\tau = O(1)$, manifesting nothing reflecting the solvability argument leading to (\ref{eq:sigma1solv}).

Notwithstanding their deficiencies, the simple approximations 
of sections~\ref{sect:tFourexp} and~\ref{sect:Texp} are instructive both for familiar reasons (namely for the transparency and simplicity of their analysis) and because they can immediately be analytically continued to assess their relevance to the evolution of singularities in the complex plane; see Figure~\ref{fig:YvsT} (in which the dashed line is derived from the approximation of section~\ref{sect:tFourexp}).



\subsection{Complex-singularity dynamics}\label{sec:comp}

The motion of the singularity necessarily breaks down into the three time scales of Appendix~\ref{sec:appgenflat} and here we adopt the specific initial condition $v = \alpha - \epsilon\cos x$. We set $x = iy$ and describe the location of the nearest singularity (corresponding to $v = 0$) on the positive imaginary axis.

On the first time scale ($t = O(1)$), we set $y = \log(1/\epsilon) + Y$ in (\ref{eq:v}) and have the full balance
\begin{equation}
    v\frac{\partial v}{\partial t} = -\left( v \frac{\partial^2 v}{\partial Y^2} - 2 \left(\frac{\partial v}{\partial Y}\right)^2  \right) - v, \label{eq:backdiff}
\end{equation}
now to be solved as an initial value problem subject to
\begin{eqnarray}
&&\text{at } t = 0, \qquad \qquad     v = \alpha - \frac{1}{2}e^{Y},  \\
&&\text{as } Y \to -\infty, \qquad   v \sim \alpha - t - \frac{1}{2}e^{Y-t}. \label{eq:backdiffcond2}
\end{eqnarray}
This not analytically solvable (and in effect a more difficult problem than the original PDE with nearly flat initial data), but does establish that to leading order the singularity location is fixed: $y \sim \log(1/\epsilon)$. This problem can, however, be solved analytically in the limits $t \to 0^{+}$ and $t \to \alpha^{-}$; using the techniques of~\cite{FKW}, we find
\begin{eqnarray}
y &\sim & \log(2\alpha/\epsilon) + \left( 2t \log(1/t)\right)^{1/2}, \qquad t \to 0^{+}, \label{eq:y1stscalel1} \\
y &\sim & \log(2/\epsilon) + \alpha + \log(\alpha - t), \quad\qquad\hspace{0.1cm} t \to \alpha^{-}, \label{eq:y1stscalel2}
\end{eqnarray}
so the singularity moves away from the real axis at early times (with unbounded speed as $t \to 0^{+}$) before reversing.

The results on the second and final time scales follow immediately from the above real line results. On the second time scale ($t = t_c + \epsilon T$, $T = O(1)$), since (see (\ref{eq:w0s}) and (\ref{eq:phi12Four})--(\ref{eq:Fourbetagamma}))
\begin{equation*}
    w_0 = -T - e^{-\alpha}\left(\cosh y - 1   \right)
\end{equation*}
the singularity location satisfies
\begin{equation}
    \cosh y \sim 1 + e^{\alpha}(-T)  \label{eq:y2ndscaleeq}
\end{equation}
so that
\begin{equation}
    y \sim \log\left( 1 + e^{\alpha}(-T) + \left( 2e^{\alpha}(-T) + e^{2\alpha}(-T)^2  \right)^{1/2}   \right).  \label{eq:y2ndscale}
\end{equation}
If follows from (\ref{eq:y2ndscale}) that
\begin{eqnarray}
&&    y  \sim  \log(-2T) + \alpha, \qquad T \to -\infty,  \label{eq:y2ndscalel1} \\
&&    y  \sim  \left(2 e^{\alpha}(-T)\right)^{1/2}, \qquad T \to 0^{-}, \label{eq:y2ndscalel2}
\end{eqnarray}
where (\ref{eq:y2ndscalel1}) exhibits the necessary matching with (\ref{eq:y1stscalel2}).

On the final time scale ($\tau = - \epsilon\log(-T)$, $\tau = O(1)$), since (see (\ref{eq:G0G1s}))
\begin{equation*}
    G_0 = 1 - \sigma_1 y^2/(-T), \qquad \sigma_1 = \frac{\gamma_1}{1 + 8 \gamma_1 \tau},
\end{equation*}
with $\gamma_1 = 1/(2e^{\alpha})$ the singularity location satisfies
\begin{equation}
    y \sim \left(  2e^{\alpha}(-T)\right)^{1/2}\left(1 - 4 \epsilon e^{-\alpha}\log(-T)  \right)^{1/2},  \label{eq:y3rdscale}
\end{equation}
which matches with (\ref{eq:y2ndscalel2}) for $\epsilon \vert \log (-T) \vert \ll 1$, while
\begin{equation*}
  y \sim \left( 8(t_c - t)\log(1/(t_c - t))  \right)^{1/2}, \qquad t \to t_c^{-},  
\end{equation*}
accordingly describes impingement onto the real axis.

\section{Flatness of the solution on the real line and proximity of the nearest singularity}\label{sec:flatness}

Note from Figure~\ref{fig:schematic} that the solutions to the $u$-equation (\ref{eq:u}) that we consider attain their maximum at $x = 0$ and their minimum at $x = \pm \pi$. Therefore, as an indication of the `flatness' of the solution profile, we consider the quantity $u(0,t) - u(\pi,t) := f(t)$, which is the relative height of the peak of the solution on $[-\pi, \pi]$. We shall consider the relation between the flatness of the solution profile and the distance of the singularities to the real axis. 

Let the $a_k(t)$ denote the Fourier coefficients of the solution in the $u$ variable. Since the solution is real (prior to blow up) and even, $a_k = a_{-k}$ and thus the flatness of the solution is given by
\begin{align}
f(t) = \sum_{k = -\infty}^{\infty} a_k(t) - \sum_{k = -\infty}^{\infty}(-1)^ka_k(t) 
= 4\sum_{k=0}^{\infty} a_{2k+1}(t).   \label{eq:flatcoeffs}
\end{align}
If $f'<0$, the solution becomes more flat while if $f'>0$ the solution becomes steeper on $x \in [-\pi, \pi]$. 
Using the residue theorem, it follows that
\begin{eqnarray}
a_k(0) &=& \frac{1}{2\pi}\int_{-\pi}^{\pi}\frac{e^{-ikx}}{\alpha - \epsilon\cos x}dx\notag  \\
&=& \frac{\left( \frac{\alpha}{\epsilon} + \sqrt{\left(\frac{\alpha}{\epsilon}\right)^2 - 1} \right)^{-\vert k \vert}}{\epsilon\sqrt{\left(\frac{\alpha}{\epsilon}\right)^2 - 1}} = \frac{1}{\alpha}\left(\frac{\epsilon}{2\alpha}\right)^{\vert k \vert} + O(\epsilon^{\vert k \vert +2}), \qquad \epsilon \to 0.  \label{eq:ICcoeffs}
\end{eqnarray}
The coefficients satisfy
\begin{equation}
    a_k' = -k^2 a_k + b_k,  \label{eq:acoeffseq}
\end{equation}
where $b_k$ is the $k$-th Fourier coefficient of $u^2$. 

\begin{figure}[htb]  
\centering
	\mbox{	
	\includegraphics[width=0.495\textwidth]{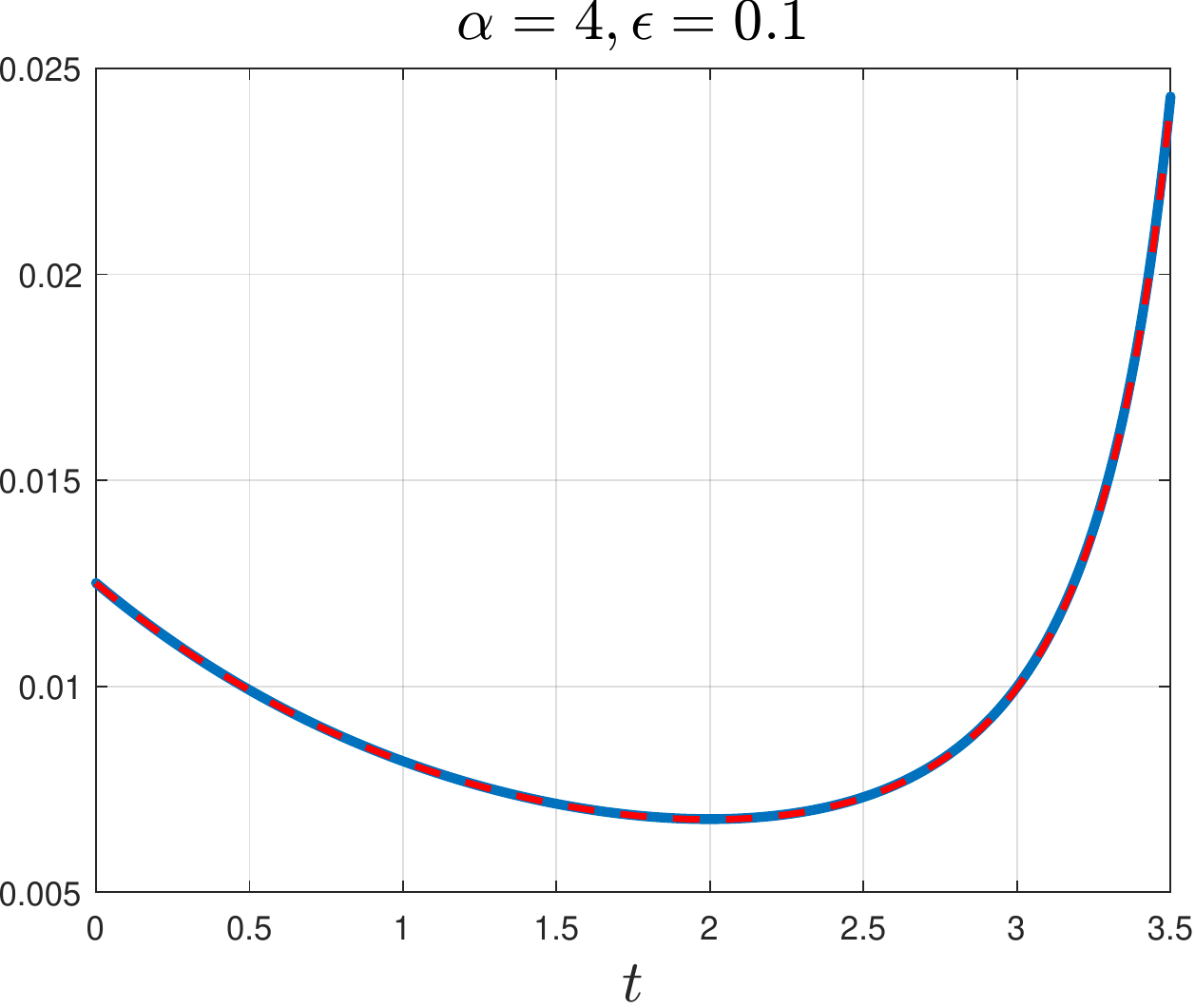}
	\includegraphics[width=0.495\textwidth]{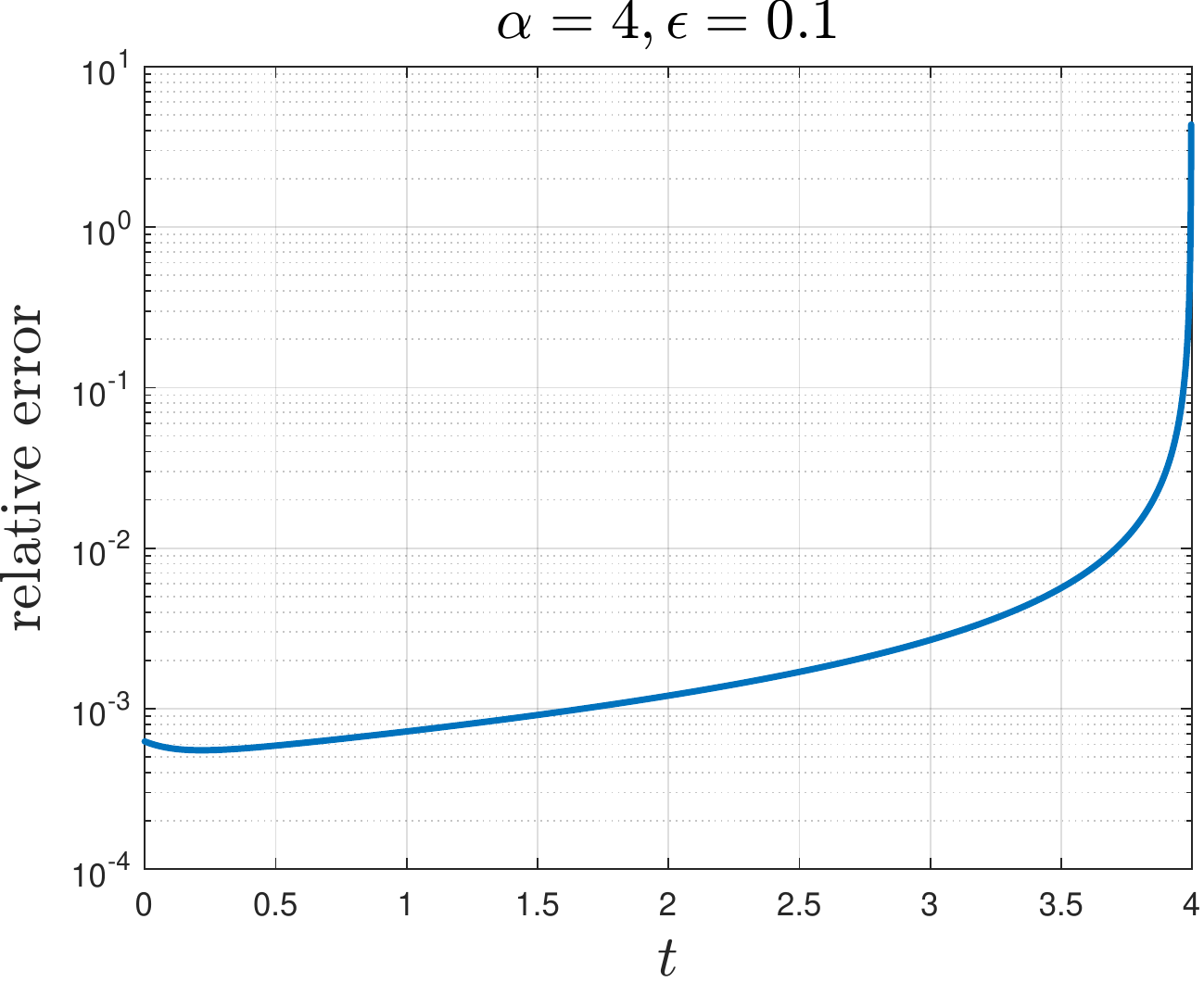}}
	\caption{  Left: The numerically computed flatness $f(t) = u(0,t)-u(\pi,t)$ (solid blue curve) compared to the approximation $f(t) \sim 4a_1$, with $a_1$ given in (\ref{eq:a0a1flat}) (red dashed curve). Right: the relative error of the approximation $f(t) \sim 4a_1$, which, as expected, ceases to be valid as blow up is approached.
	}  \label{fig:flatness}
\end{figure}

It follows from (\ref{eq:ICcoeffs})--(\ref{eq:acoeffseq}) that, to leading order, as $\epsilon \to 0$, and away from the blow-up time, 
\begin{equation}
    a_0(t) \sim \frac{1}{\alpha - t}, \qquad a_1(t) \sim \frac{\epsilon e^{-t}}{2(\alpha - t)^2}, \label{eq:a0a1flat}  
\end{equation}
therefore
\begin{equation*}
    f'(t) \sim 4a_1'
    \sim \frac{2\epsilon \,e^{-t}(t - (\alpha - 2))}{(\alpha - t)^3}.
\end{equation*}
Hence, on the real line, the solution switches from flattening to steepening at  $t \sim \alpha - 2$ if $\alpha > 2$ and the minimal flatness is $f(\alpha-2)\sim 4a_1(\alpha - 2) \sim \epsilon\, e^{2-\alpha}/2$. If $\alpha < 2$, the solution does not flatten at all but steepens from $t=0$ until the blow-up time. Figure~\ref{fig:flatness} confirms the validity of the flatness approximation $f(t) \sim 4a_1(t)$ away from the blow-up time.

In the complex plane, we deduce from the asymptotic approximation (\ref{eq:y1stscalel1}) that, regardless of the parameter values (provided $\alpha \ll \epsilon$), the singularities initially move away from the real axis and turn around at $t \sim \exp(-1) \approx 0.37$, which is consistent with the numerical results in Figure~\ref{fig:YvsT}. 
This illustrates that there is not a simple correspondence between the distance of the singularity from the real axis and the flatness of the solution. 

Indeed, it follows from (\ref{eq:flatcoeffs}) and the rapid decay of the Fourier coefficients (away from the blow-up time) that the flatness of the solution is determined by the behaviour of the low-order Fourier coefficient $a_1$. On the other hand,  the distance of the singularity from the real axis is determined by the behaviour of the high-order Fourier modes. Hence, the discrepancy between the flatness of the solution and the proximity of the singularity is due to the qualitatively different evolution of the lower and higher order Fourier modes.

For the high-order modes, the diffusion term $-k^2a_k$ dominates the nonlinear term $b_k$ for small $t$ (again, due to the rapid decay of the Fourier coefficients). Therefore $a_k' < 0$, the high-order Fourier coefficients decay for small $t$ and therefore the singularity initially moves away from the real axis, regardless of the parameter values satisfying $\epsilon \ll \alpha$. The high-order Fourier coefficients start increasing as the singularity turns around, which occurs at $t \sim 0.37$, according to (\ref{eq:y1stscalel1}).

\end{document}